\documentclass[12pt]{article}
\usepackage{amsmath,amssymb}
\usepackage{verbatim}
\usepackage{dsfont,bm}
\usepackage{color}
\usepackage{graphicx}
\usepackage{amsthm}
\usepackage{enumerate}

\usepackage{xcolor}
\definecolor{darkgreen}{rgb}{0,0.75,0}
\definecolor{darkred}{rgb}{0.75,0,0}
\definecolor{darkmagenta}{rgb}{0.5,0,0.5}
\usepackage[colorlinks,citecolor=darkgreen,linkcolor=darkred,urlcolor=darkmagenta]{hyperref}
\newtheorem{theorem}{Theorem}[section]

\newtheorem{lemma}[theorem]{Lemma}
\newtheorem{lem}[theorem]{Lemma}

\newtheorem{prop}[theorem]{Proposition}
\newtheorem{definition}[theorem]{Definition}

\newtheorem{remark}[theorem]{Remark}
\newtheorem{example}[theorem]{Example}

\numberwithin{equation}{section}

\def\be{\begin{equation}}
\def\ee{\end{equation}}
\def\bes{\begin{equation*}}
\def\ees{\end{equation*}}

\newcommand{\Mod}[0]{\operatorname{Mod}}
\newcommand{\mr}[1]{{\tt \href{http://www.ams.org/mathscinet-getitem?mr=#1}{MR#1}}}
\newcommand{\arxiv}[1]{{\tt \href{http://arxiv.org/abs/#1}{arXiv:#1}}}
\newcommand{\set}[1]{\left\{ #1 \right\}}

\newcommand{\Sett}[2]{\left\{ #1  : \, #2 \right\}}

\newcommand{\abs}[1]{{\left\vert\kern-0.25ex #1
		\kern-0.25ex\right\vert}}
\newcommand\norm[1]{\left\lVert#1\right\rVert} 
\newcommand{\one}{\mathds{1}} 
\newcommand{\loc}[0]{\operatorname{loc}}

\newcommand{\compl}{c}

\def\sA {{\mathcal A}}  \def\sC {{\mathcal C}}
\def\sD {{\mathcal D}} \def\sE {{\mathcal E}} \def\sF {{\mathcal F}}
 \def\sH {{\mathcal H}} 
  
\def\sM {{\mathcal M}}  
  
\def\sS {{\mathcal S}}  
 \def\sW {{\mathcal W}} \def\sX {{\mathcal X}}

 \def\bE {{\mathbb E}} 
\def\bG {{\mathbb G}}  
  
 \def\bN {{\mathbb N}} 
\def\bP {{\mathbb P}}  \def\bR {{\mathbb R}}
\def\bS {{\mathbb S}}  \def\bU {{\mathbb U}}
  \def\bX {{\mathbb X}}

\def\med{\medbreak\noindent}

\def\sm{\smallskip\noindent}

\def\ignore#1{}
\def\noi{\noindent}

\def\ol{\overline}           

\def\Gam{\Gamma}

\def\to {\rightarrow}

\def\dint{\int\kern-.6em\int}

\newcommand\restr[2]{{
		\left.\kern-\nulldelimiterspace 
		#1 
		\vphantom{\big|} 
		\right|_{#2} 
	}} 
	
	\def\diam{{\mathop{{\rm diam }}}}
	\def\dist{{\mathop {{\rm dist}}}}

	\def\Cap{\operatorname{Cap}}
	\def\Mod{{\mathop{{\rm Mod}}}}

	\newcommand{\on}[1]{\operatorname{ #1}}
	
	\def\wt{\widetilde}
	\def\wh{\widehat}
	
	\def\be{\begin{equation}}
	\def\ee{\end{equation}}
	\def\bes{\begin{equation*}}
	\def\ees{\end{equation*}}
	\def\ba{\begin{align}}
	\def\ea{\end{align}}
	\def\xxea{\end{align}}
\def\bas{\begin{align*}}
\def\eas{\end{align*}}

\def\proof{{\smallskip\noindent {\em Proof. }}}
\def\qed{{\hfill $\square$ \bigskip}}

\definecolor{dgreen}{rgb}{0, 0.6, 0.1}
\definecolor{dblue}{rgb}{0, 0.0, 0.6}
\definecolor{vdblue}{rgb}{0,.08, 0.45}
\definecolor{dred}{rgb}{0.7, 0.0, 0.0}
\definecolor{vdblue}{rgb}{0,.08, 0.45}

\definecolor{purple}{rgb}{0.6, 0.0, 0.6}
\definecolor{mytext}{rgb}{0.1, 0.1, 0.1}

\begin{document}
	
	\font\titlefont=cmbx14 scaled\magstep1
	\title{\titlefont  Quasisymmetric uniformization and  heat kernel estimates}  
	\author{
		Mathav Murugan\footnote{Research partially supported by NSERC (Canada) and the Pacific Institute for the Mathematical Sciences} 
	}
	\maketitle
	\vspace{-0.5cm}
	{\centering Dedicated to Professor Laurent Saloff-Coste on the occasion of his 60th birthday.\par}
	
\begin{abstract}
We show that the circle packing embedding  in $\bR^2$ of a one-ended, planar triangulation with polynomial growth is quasisymmetric  if and only if the simple random walk on the graph satisfies sub-Gaussian heat kernel estimate with spectral dimension two.
Our main results provide a new family of graphs and fractals that satisfy sub-Gaussian estimates and Harnack inequalities. 

\vskip.2cm
\noindent {\it Keywords:} 
Quasisymmetry, Uniformization, Circle packing, Sub-Gaussian estimate,  Harnack inequality.
\vskip.2cm \noindent {\it Subject Classification (MSC 2010): 60J45, 51F99 }
\end{abstract}

\section{Introduction}

The classical uniformization theorem implies that  a Riemann surface that is homeomorphic to the $2$-sphere is conformally equivalent to $\bS^2$. 
Therefore, the Brownian motion associated with a conformal metric on such a Riemann surface can be viewed as a time change of the Brownian motion on $\bS^2$.
In this work, we show that a similar property holds for Brownian motion on metric spaces  with a notion of  `generalized conformal map' to $\bS^2$. 
Furthermore, as we shall see, the availability of such a generalized conformal map allows us to obtain Harnack inequalities and heat kernel estimates  for diffusions and random walks. This work explores a new relationship between quasiconformal geometry of fractals and diffusion on fractals.

Quasisymmetric maps are a fruitful generalization of conformal maps.  Quasisymmetric maps  were introduced by Beurling and Ahlfors, and were studied as boundary values of quasiconformal self maps of the upper half plane \cite{BA}. Heinonen's book \cite{Hei} is an excellent reference on quasisymmetric maps. We recall the definition due to \cite{TV} below.
\begin{definition}\label{d:qs}
	{\rm	 A \emph{distortion function} is a homeomorphism of $[0,\infty)
		$ onto itself. 
		Let $\eta$ be a distortion function. 
		A map $f:(\sX_1,d_1) \to (\sX_2,d_2)$ between metric spaces is said to be
		\emph{$\eta$-quasisymmetric}, if $f$ is a homeomorphism and
		\[
		\frac{d_2(f(x),f(a))}{d_2(f(x),f(b))} \le \eta\left(\frac{d_1(x,a)}{d_1(x,b)}\right)  
		\]
		for all triples of points $x,a,b \in \sX_1, x \neq b$. We say $f$ is a \emph{quasisymmetry} if it is $\eta$-quasisymmetric for some distortion function $\eta$.
		We say that metric spaces $(\sX_1,d_1)$ and $(\sX_2,d_2)$ are quasisymmetric, if there exist
		s a quasisymmetry $f:(\sX_1,d_1) \to 
		(\sX_2,d_2)$.
		We say that  metrics $d_1$ and $d_2$  on $\sX$ are \emph{quasisymmetric}, if the identity map $\operatorname{Id}:(\sX,d_1) \to (\sX,d_2)$ is a quasisymmetry.
		We say that a (not necessarily onto) map $f:(\sX_1,d_1) \to (\sX_2,d_2)$ is a quasisymmetric embedding if $f: (\sX_1,d_1) \to (f(\sX_1),d_2)$ is a quasisymmetry.
		 }
\end{definition}

An important motivation behind generalizations of conformal structures arise from geometry of hyperbolic spaces. We refer the reader to  \cite[Chapters 1 and 9]{GMP}, 
and the ICM surveys of Bonk \cite{Bon} and of Kleiner \cite{Kle} for a good exposition of the relationship between quasisymmetric maps, uniformization and the geometry of hyperbolic spaces. A fundamental relationship between hyperbolic spaces and quasisymmetry is the following: two Gromov hyperbolic spaces \cite{GH} are quasi-isometric if and only if their boundaries are quasisymmetric. This observation arises from Mostow's celebrated work on rigidity theorem -- see \cite{Pau,BSc} for modern formulations.

Motivated by the above considerations, one is led to the {\bf quasisymmetric uniformization problem}: ``What conditions should be imposed on a metric space  $(\sX,d)$ so that  it is quasisymmetric to a model space $(\sM,d_\sM)$?'' There is a simple answer to the quasisymmetric uniformization problem when the model space is $\bR$ or $\bS^1$ due to Tukia and V\"ais\"al\"a \cite{TV}, \cite[Theorem 15.3]{Hei}.
We refer the reader to  \cite{BK,Raj} for important results on the quasisymmetric uniformization problem when the model space is $\bS^2$. A different combinatorial approach to equip spaces with generalized conformal structures was developed by Cannon \cite{Can}.

The primary message of our work is that  quasisymmetric uniformization for $\bR^2$ and $\bS^2$ is closely related to random walks and diffusions.
From a probabilistic viewpoint, the existence of a quasisymmetric map from a metric space to a well understood model space allows us to study diffusions and random walks.
Many properties that are relevant to random walks and diffusions can be transferred using a quasisymmetry; for example, the elliptic Harnack inequality, Poincar\'e inequality and resistance estimates \cite[Section 5]{BM1}.
More generally, changing the metric of a space  is an useful tool to study diffusions and random walks \cite{ABGN, BeSc1,BeSc2, Geo, GN,  Kig12, Lee1,Lee2}.

We denote the  graph distance on a connected graph $\bG=(V_\bG,E_\bG)$ by $d_\bG:V_\bG \times V_\bG \to [0,\infty)$. The open ball with center $x \in V_\bG$ with radius $r$ in the $d_\bG$ metric is denoted by $B_\bG(x,r):= \set{y \in V_\bG: d_\bG(x,y) < r}$.
We say that a graph $\bG=(V_\bG,E_\bG)$ is of \emph{polynomial growth} with  \emph{volume growth exponent} $d$ if there exists $C>1$ such that for any ball $B_\bG(x,r)$, its cardinality $\abs{B_\bG(x,r)}$ satisfies the estimate $C^{-1} r^d \le \abs{B_\bG(x,r)} \le C r^d$ for all $x \in V_\bG, r \ge 1$. 

Many regular fractals and their graph analogues satisfy sub-Gaussian transition probability estimates. Such estimates were first obtained for the Sierpinski gasket \cite{BP}. We refer the reader to  \cite{Bar1} for an introduction to sub-Gaussian estimates and diffusion on fractals. 
\begin{definition}
	{\rm
	We say that  a graph $\bG=(V_\bG,E_\bG)$ of {polynomial growth} with  {volume growth exponent} $d$ satisfies the \emph{sub-Gaussian estimate with walk dimension} $d_w$, if  there exists $C>1$ such that the simple random walk $(Y_n)_{n \ge 0}$ admits the following heat kernel upper and lower bounds:
	\be \label{e:usg}
	\bP^x(Y_n=y) \le \frac{C}{n^{d/d_w}} \exp \left( - \left( \frac{d_\bG(x,y)^{d_w}}{Cn}\right)^{1/(d_w-1)} \right) \hspace{2mm} \mbox{for all $x,y \in V_\bG$ and $n \ge 1$,}
	\ee 
	and
	\be \label{e:lsg}
	\bP^x(Y_n=y) + \bP^x(Y_{n+1}=y) \ge \frac{C^{-1}}{n^{d/d_w}} \exp \left( - \left( \frac{Cd_\bG(x,y)^{d_w}}{n}\right)^{1/(d_w-1)} \right), 
	\ee  
	for all $x,y \in V_\bG$ and $n \ge 1 \vee d_\bG(x,y)$, where $\bP^x$ denote the probability conditioned on the event that the random walk starts at $Y_0=x$.
	
The quantity $d_s= 2 d/ d_w$ is called the \emph{spectral dimension}. 
	}
\end{definition}
\begin{remark}
	{\rm 
 The sub-Gaussian estimates \eqref{e:usg}, \eqref{e:lsg} can be understood better by recalling some well-known consequences.
	\begin{enumerate}
		\item 
	The spectral dimension gives return probability estimates as $ \bP^x(Y_n = x) +  \bP^x(Y_{n+1} = x) \asymp n^{-d_s/2}$. In particular, the simple random walk is transient if and only if $d_s > 2$.
	\item
	These heat kernel bounds imply that the expected distance travelled by the  random walk satisfies the estimate  $\bE (d_\bG(Y_0,Y_n)) \asymp n^{1/d_w}$, and the expected time to exit a ball of radius $r$ satisfies the estimate 
	$\bE^x (\tau_{B_\bG(x,r)}) \asymp r^{d_w}$, where $\tau_{B_\bG(x,r)}$ denotes the exit time of the random walk from the ball $B_\bG(x,r)$.
	\item  The bounds $2 \le d_w \le 1+ d$ always hold \cite{Bar}.  Taking  $d_w=2$ corresponds to the classical Gaussian estimates.
		\end{enumerate}
}
\end{remark}

Our main result (Theorem \ref{t:main1}) is not stated in the introduction because its formulation requires some preparation. 
Instead, we state a consequence of the main result that relates circle packing to random walks.
Recall that
a \emph{circle packing} of a planar graph $\bG=(V_\bG,E_\bG)$ is a set of of circles
with disjoint interiors $\{C_v\}_{v\in V_\bG}$ in the plane $\bR^2$ such that two 
circles are
tangent if and only if the corresponding vertices form an edge.  This provides an 
embedding $f_{\on{CP}}: V_\bG \to \bR^2$ which sends the vertices to the centres of the corresponding circles, and induces a  \emph{circle packing metric} $d_{\on{CP}} : V_\bG \times V_\bG \to [0,\infty)$ defined by $d_{\on{CP}}(x,y) = \abs{ f_{\on{CP}}(x) - f_{\on{CP}}(y)}$, where $\abs{a}$ denotes the Euclidean norm of $a \in \bR^2$.

\begin{theorem} \label{t:main2}
Let $\bG=(V_\bG,E_\bG)$ be one-ended, planar triangulation of polynomial growth with volume growth exponent $d$. Then the following are equivalent
	\begin{itemize}
		\item [(a)] The circle packing metric $d_{\on{CP}}$ and the graph metric $d_{\on{\bG}}$ are quasisymmetric.
		\item [(b)] The simple random on $\bG$ satisfies sub-Gaussian estimate with walk dimension $d_w=d$ (or equivalently, the spectral dimension $d_s=2$).
	\end{itemize}
\end{theorem} 
\begin{remark}{\rm 
Although Theorem \ref{t:main2} only applies for triangulations, it is more widely applicable due the reasons outlined below.
\begin{enumerate}[(a)]
\item Consider a one-ended planar graph $\bG$, such that both $\bG$ and its planar dual $\bG^\dagger$ are of bounded degree. Then the circle packing embedding in (a) of Theorem \ref{t:main2} can be replaced by a more general notion of good embedding (see Lemma \ref{l:good} and Theorem \ref{t:main1}).
This extends Theorem \ref{t:main2} to a larger family of graphs, which for example contains quadrangulations.
\item Alternately, one could use the face barycenter triangulation of \cite{CFP2} to obtain a triangulation $\wt{\bG}$ such that $\wt{\bG}$ is quasi-isometric to $\bG$. The stability results for sub-Gaussian estimates \cite{BB} allow us to see that sub-Gaussian estimates for $\bG$ are equivalent to sub-Gaussian heat kernel estimates for the triangulation $\wt{\bG}$.


\item As shown in the work of Bonk and Kleiner,  planar surfaces can be triangulated and properties of the metric space can be transferred between the discrete triangulations and  metric surfaces \cite[Theorem 11.1 and Section 8]{BK}.  Furthermore, \cite{BK} use circle packing to construct quasisymmetric maps on metric $2$-spheres as a limit of circle packings of finer and finer triangulations of the metric space. 
A key ingredient in \cite{BK} is that the circle packing embeddings are uniformly quasisymmetric.
This is one of the motivations behind studying quasisymmetry of the circle packing embedding on graphs.
\end{enumerate}

}\end{remark}

We mention some further motivations behind this work and provide some context.
The Uniform Infinite Planar Triangulation (UIPT) and $\sqrt{8/3}$-Liouville quantum gravity are conjectured to have spectral dimension $2$. This conjecture should be interpreted in a weaker sense as sub-Gaussian estimates are unrealistically strong as these space do not satisfy the volume doubling property. Although satisfactory heat kernel bounds are still unknown, there has been spectacular recent progress in showing that UIPT and several other random planar maps have spectral dimension $2$ by obtaining bounds on resistances and exit times \cite{GM, GHu}.  We expect that some of the methods developed here will be useful to obtain heat kernel bounds for random planar maps.
The proof of Theorem \ref{t:main2} solves a special case of the `resistance conjecture'  when the graph is a one-ended, planar triangulation with polynomial growth  and has spectral dimension $d_s=2$; see Remark \ref{r:main}(ii). 
Therefore this work can be considered as evidence towards the resistance conjecture when $d_s=2$.
 We recall that the case $d_s<2$ of the resistance conjecture has been solved in \cite{BCK}.

Our main result and examples provide  non-trivial examples of spaces with conformal walk dimension two -- see \cite[Remark  5.16(3)]{BM1} for the definition of conformal walk dimension. The only previously known (to the best of the author's knowledge)  non-trivial example of a space with conformal walk dimension two is the Sierpinski gasket,  which follows from Kigami's work on the `harmonic Sierpinski gasket' \cite{Kig08}. It is not known if the conformal walk dimension of the standard Sierpinski carpet is two -- see  \cite[Section 9]{Ka} for related questions.

The outline of this work is as follows. 
In Section \ref{s:prelim}, we recall some some background on modulus of curve families, Dirichlet forms, cable processes, and circle packing of planar graphs. 
The proof of the implication `$\operatorname{(b) \implies (a)}$' in Theorem \ref{t:main2} follows an old approach of Heinonen and Koskela \cite{HK}. The main ingredients required to carry out this approach are the annular quasi-convexity for the graph metric $d_\bG$, and the Loewner property for the circle packing metric $d_{\on{CP}}$.
In Section \ref{s:aqc}, we introduce the notion of annular quasi-convexity at large scales, and obtain a sufficient condition using the Poincar\'e inequality and capacity bounds. In Section \ref{s:loew}, we recall the definition of the Loewner property, and obtain  the Loewner property for the circle packing embedding under mild conditions. The key tool to show the Loewner property is the Poincar\'{e} inequality for the circle packing embedding obtained in \cite{ABGN}.
 In Section \ref{s:qs}, we  obtain quasisymmetry between the graph and circle packing metrics using above mentioned annular quasi-convexity at large scales, and the Loewner property.  In Section \ref{s:hk}, we show that quasi-symmetry of the circle packing embedding implies sub-Gaussian bounds on the heat kernel using the methods of \cite{BM1}.
Finally, as an application of the main results, we obtain sub-Gaussian estimates for a new family of graphs in Section \ref{s:xm}. Furthermore, we show that these methods can also be applied to obtain heat kernel bounds for diffusions on a family of fractals that are homeomorphic to $\bS^2$. 
\section{Preliminaries} \label{s:prelim}
We recall some basics on Dirichlet forms \cite{FOT,CF}, modulus of curve families \cite{Hei,HKST}, and cable systems and processes \cite{BBI, Fol}.

\subsection{Upper gradient}
Let $(\sX,d)$ be a complete metric space.
We recall the notion of rectifiability and line integration. 
By a \emph{curve} we mean  a continuous map $\gamma:I \to \sX$ of an interval $I \subset \bR$ into $\sX$.
A \emph{subcurve} of $\gamma$ is the restriction $\restr{\gamma}{J}$ of $\gamma$ to a subinterval $J \subset I$.
 We sometimes abuse notation and abbreviate the image $\gamma(I)$ by $\gamma$. 
 If $I=[a,b]$, then the length of the curve $\gamma:I \to \sX$ is 
 \[
 L(\gamma)= \on {len
gth}(\gamma) = \sup \sum_{i=1}^n d(\gamma(t_i),\gamma(t_{i+1})),
 \]
where the supremum is over all finite sequences $a=t_1 \le t_2 \le \ldots \le t_n \le t_{n+1}=b$. If $I$ is not compact, then we set
\[
L(\gamma) = \sup_{J} L \left( \restr{\gamma}{J}\right),
\]
where $J$ varies over all compact subintervals of $I$.
 We say $\gamma$ is \emph{rectifiable} if $L(\gamma) < \infty$. 
Similarly, a curve $\gamma:I \to \sX$ is \emph{locally rectifiable} if its restriction to each 
compact subinterval of $I$ is rectifiable. 

Any rectifiable curve admits a unique extension $\ol{\gamma}: \ol{I}  \to \sX$. If $I$ is unbounded the extension is understood in a generalized sense. From now on, given any rectifiable curve $\gamma$ we automatically consider its extension $\ol{\gamma}$ and do not distinguish in notation. 
Any rectifiable curve admits a natural \emph{arc length parametrization}  defined as the unique $1$-Lipschitz map  $\gamma_s:[0,L(\gamma)] \to \sX$ such that $\gamma= \gamma_s \circ s_\gamma$ where $s_\gamma:[a,b] \to [0,L(\gamma)]$ is the length function $s_\gamma(c)= L( \restr{\gamma}{[a,c]})$. 
If $\gamma$ is a rectifiable curve in $\sX$, the \emph
{line integral} over $\gamma$ of each non-negative Borel function $\varrho:\sX \to [0
,\infty]$ is 
\[
\int_{\gamma} \varrho\,ds =
 \int_0^{L(\gamma)} \varrho \circ \gamma_s(t) \, dt.
\]

We recall the notion of upper gradient \cite[p. 152]{HKST}.
\begin{definition}[Upper gradient] \label{d:ug}
{\rm Let $(\sX,d)$ be a metric space and let $u:\sX \to \bR$ be a function.
A Borel function $
\rho: \sX \to \bR$ is said to be 
an \emph{upper gradient} of $u$, if
\[
\abs{u(\gamma(a)) - 
 u(\gamma(b))} \le \int_
{\gamma} \rho \, ds,
\]
for every rectifiable curve $\gamma:[a,b] \to \sX$.}
\end{definition}
\subsection{Modulus of a curve family}
Let $(\sX,d,\mu)$ a complete, metric measure 
space such that $\mu$ is a Radon measure with full support.
Let $\Gamma$ be a family of curves in $\sX$ and $p>0$. The \emph{$p$-modulus} of $\Gamma$ is defined as 
\[
\Mod_p(\Gam) = \inf_\varrho \int_{\sX} \varrho^p \, d\mu,
\]
where the infimum is taken over all nonnegative Borel functions $\varrho:\sX \to [0,\infty]$ satisfying 
\be \label{e:admiss}
\int_{\gamma} \varrho \, ds \ge 1
\ee
for all locally rectifiable curves $\gamma \in \Gamma$. Functions satisfying \eqref{e:admiss} are called \emph{admissible functions} or admissible metrics  for $\Gam$.
We recall the basic properties of  the $p$-modulus  \cite[p. 128]{HKST}. They are
\begin{align}
\label{e:bpm1}\Mod_p (\emptyset) & =
 0, \\
\label{e:bpm2}\Mod_p (\Gamma_1) & \le \Mod_p(\Gamma_2) \hspace{2mm} \mbox{for all $\Gamma_1 \subset \Gamma_2$,} \\
\label{e:bpm3}\Mod_p\left( \cup_{i=1}^\infty \Gamma_i \right) & \le \sum_{i=1}^\infty \Mod_p(\Gam_i),
\\
\label{e:bpm4}\Mod_p(\Gamma) & \le \Mod_p(\Gam_0),
\end{align}
whenever $\Gamma_0$ and $\Gamma$ are two curve families such that each curve $\gamma \in \Gamma$ has a subcurve $\gamma_o \in \Gamma_0$, \emph{i.e.} $\Gamma$ has fewer and longer curves than $\Gamma_0$.

 We will exclusively use $2$-modulus and therefore will abbreviate $\Mod_2(\cdot)$ by $\Mod(\cdot)$. Similarly, by modulus we mean $2$-modulus.
By $\Mod(E,F;U)$ we denote the modulus of the family of all curves in a subset $U$ of $\sX$  joining two disjoint subsets $E$ and $F$ of $U$.  We abbreviate $\Mod(E,F;\sX)$ by $\Mod(E,F)$.

\subsection{Dirichlet forms}
We recall some standard notions concerning Dirichlet forms and refer the reader to \cite{FOT,CF} for a detailed exposition. Let $(\sX,d,\mu)$ be a locally compact, separable, metric measure space, where $\mu$ is a Radon measure with full support.
Let $(\sE,\sF)$ be a strongly local, regular Dirichlet form
 on $L^2(\sX,m)$ -- see \cite[Sec. 1.1]{FOT}. Associated with this 
form  $(\sE,\sF)$, there exists an $\mu$-symmetric Hunt process  $\bX= (\Omega, \sF_\infty,\sF_t,X_t, \bP_x)$ \cite[Theorem 7.2.1]{FOT}.
We denote the \emph{extended Dirichlet space} by $\sF_e$ \cite[Theorem 1.5.2]{FOT}.
Recall that the 
Dirichlet form $(\sE,\sF)$ is \emph{recurrent} if and only if $1 \in \sF_e$ and $\sE(1,1
)=0$ \cite[Theorem 1.6.3]{FOT}. For $f \in \sC_c(\sX) \cap \sF$,  the \emph{energy measure} is defined as the unique Borel measure $d\Gamma(f,f)$ on $\sX$ that satisfies
$$ \int g d\Gam(f,f) = 2 \sE(f,fg) - \sE(f^2,g), \hspace{3mm}\mbox{ for all $g \in \sF \cap \sC_c(\sX)$.} $$
This notion can be extended to all functions in $\sF_e$ and
we have
$$ \sE(f,f) = \int_{\sX} d\Gam(f,f). $$
This follows from \cite[Lemma 3.2.3]{FOT} with a caveat that our definition of $\Gamma(f,f)$ is different from \cite{FOT} by a factor $1/2$.

Let $(\sX,d)$ be a metric space equipped  with a  strongly local Dirichlet form $(\sE,\sF)$ on $L^2(\sX,\mu)$.
We call  $(\sX,d,\mu,\sE,\sF)$ a {\em  metric  measure space with Dirichlet form}, or {\em MMD space}.
We define \emph{capacities} for a MMD space $(\sX,d,\mu,\sE,\sF)$  as follows. For 
a non-empty open subset $D \subset \sX$, let $\sC_c(D)$ denote the space of all continuous functions with compact support in $D$. Let $\sF_D$ denote the closure of $\sF \cap \sC_c(D)$ with 
respect to the $\left(\sE(\cdot,\cdot)+ \langle \cdot,\cdot \rangle_{L^2(\mu)}\right)^{1/2}$-norm. 
By $A \Subset D$, we mean that the closure of $A$ is a compact subset of $D$.
For $A \Subset D$ we set
\be \label{e:capdef}
\Cap_D(A) = \inf\{ \sE(f,f): f \in \sF_D \mbox{ and  $f \ge 1$ in a neighbourhood of $A$} \}.
\ee
The following \emph{domain monotonicity} of capacity is clear from the definition: if $A_1 \subset A_2 
\Subset D_1 \subset D_2$ then
\be \label{e:capmon}
\Cap_{D_2}(A_1) \le \Cap_{D_1}(A_2). 
\ee


The following upper bound on the capacity will play an important role in this work. 
\begin{definition} \label{d:capub}
{\rm	We say that a MMD space satisfies \hypertarget{cap}{$(\on{cap}_{\le})$}, if
	there exists $C,M\ge 1$ such that for all $x \in \sX, r \ge 1$, we have 
	\be  \tag*{$(\on{cap}_{\le})$}
	\Cap_{B(x,Mr)} (B(x,r)) \le C.
	\ee}
\end{definition}
We record an easy consequence of \hyperlink{cap}{$(\on{cap}_{\le})$} below.
\begin{lem} \label{l:series} 
	Let $(\sX,d,\mu,\sE,\sF)$ be a MMD space that satisfies \hyperlink{cap}{$(\on{cap}_{\le})$}. 
	Then  there exists $C_1>0$ such that for all 
	$x \in \sX, r \ge 1, R \ge Mr$, we have
	\[
	\Cap_{B(x,R)} (B(x,r)) \le C_1 \left(\log (R/r)\right)^{-1},
	\]
	where $M>1$ is the constant in Definition \ref{d:capub}.
\end{lem}
\proof
Let $r \ge 1, R \ge Mr$.
Let $k \ge 1$ be the largest integer such that $R \ge M^k r$.
By \hyperlink{cap}{$(\on{cap}_{\le})$}, there exists $f_i \in \sC_c\left(B(x,M^i r) \right),  i=1,\ldots,k$ such that $f \ge 1$ in $B(x,M^{i-1},r)$ and
\be \label{e:ser1}
\sE(f_i,f_i) \le 2 \Cap_{B(x,M^i r)} \left(B(x,M^{i-1}r)\right) \le 2 C.
\ee
Since $(\sE,\sF)$ is strongly local, we have $\sE(f_i,f_j)=0$ for $i \neq j$, and therefore $f= \left(\sum_{i=1}^k f_i \right)/k$ satisfies $f \ge 1$ in $B(x,r)$, $f \in \sC_c(B(x,R))$, and
\[
\sE(f,f) =  \frac{1}{k^2} \sum_{i=1}^k \sE(f_i,f_i) \le 2C \frac{1}{k} \le C_1 \left(\log (R/r)\right)^{-1}.
\]
We use \eqref{e:ser1} in the inequality above.
\qed

One often needs lower bounds on the capacity as well. 
The following Poincar\'e inequality is used to obtain such lower bounds.
\begin{definition}\label{d:poin} {\rm
Let $\Psi:(0,\infty) \to (0,\infty)$.
We say that a MMD space $(\sX,d,\mu,\sE,\sF)$ satisfies the  \emph{Poincar\'e inequality} \hypertarget{pi}{$\on{PI}(\Psi)$} if there exist  $C,M \ge 1$ such that for all $f \in \sF, x \in \sX, 0 < r < \diam (\sX)/2$, we have
\be \tag*{$\on{PI}(\Psi)$} 
\int_{B(x,r)} \abs{f(y)-f_{B(x,r)}}^2 \, d\mu(y) \le C \Psi(r) \int_{B(x,Mr)} d\Gam(f,f),
\ee
where $f_{B(x,r)} = \frac{1}{\mu(B(x,r))} \int_{B(x,r)} f\,d\mu$.} By $\on{PI}(2)$, we mean $\on{PI}(\Psi)$ with $\Psi(r)=r^2$.
\end{definition}

We recall a useful condition to check if a function belongs to the extended Dirichlet space.
\begin{lem}\cite[Lemma 2]{Sch} \label{l:extd}
	Let $u \in L^0(\sX,\mu)$ and let $\set{u_n}_{n\in \bN}$ be a sequence in $\sF$ such that
$\lim_{n \to \infty} u_n = u$ $\mu$-almost everywhere and $\liminf_{n \to \infty} \sE(u_n,u_n) < \infty$. Then $u \in \sF_e$, and $\sE(u,u) \le \liminf_{n \to \infty} \sE(u_n,u_n)$. 
\end{lem}

We recall the notion of \emph{harmonic} functions associated to a MMD space $(\sX,d,\mu,\sE,\sF)$. We denote the local Dirichlet space corresponding to an open subset  $U$ of $\sX$ by
\[
\sF_{\loc}(U) = \Sett{u \in L^2_{\loc}(U,\mu)}{ \forall \mbox{ relatively compact open } 
V \subset U, \exists u^\# \in \sF, u = u^\# \big|_{V}  \mbox{ $\mu$-a.e.}}\]
\begin{definition}\label{d-harmonic}
	{\rm Let $U \subset \sX$ be open. A function $u: \sX \to \bR$ is \emph{harmonic} on $U$
		if  $u \in \sF_{e}$ and for any function $\phi \in \sC_c(U) \cap \sF$, we have
		\[
		\sE(u^\#, \phi)= 0.
		\]
where $u^\# \in \sF$ is such that $u^\# = u$ in the essential support of $\phi$.
	}\end{definition}	
	\begin{remark} \label{r:harm}
		{\rm  It is known that $u \in L^\infty_{\loc}(\sX,\mu)$ is harmonic in $U$ if and only if it satisfies the following property: for every relatively compact open subset $V$
				of $U$, $t \mapsto \widetilde{u}(X_{t \wedge \tau_V})$ is a uniformly integrable $\bP^x$-martingale 
				for q.e.~$x \in V$. (Here $\widetilde{u}$ is a quasi continuous version of $u$ on $V$.)
				This equivalence between the weak solution formulation in Definition \ref{d-harmonic} and the probabilistic formulation using martingales is given in \cite[Theorem 2.11]{Che}.
		}\end{remark}
It is also easy to observe that the Poincar\'e inequality extends  to functions in the local Dirichlet space $\sF_{\loc}(\sX)$. 		

We recall the notion of a heat kernel. Let  $\bX= (\Omega, \sF_\infty,\sF_t,X_t, \bP_x)$  denote the Hunt process corresponding to a MMD space $(\sX,d,\mu,\sE,\sF)$. We say that a measurable function $p: (0,\infty) \times \sX \times \sX \to [0, \infty)$ is the \emph{heat kernel} corresponding to the MMD space  $(\sX,d,\mu,\sE,\sF)$  if 
\[
\bP_x (X_t \in A) = \int_A p(t,x,y) \, \mu(dy) \hspace{2mm} \mbox{ for all $x \in \sX$ and for all Borel sets $A \subset \sX$.}
\]
We introduce a continuous time variant of the sub-Gaussian estimates in \eqref{e:usg}, \eqref{e:lsg}.
\begin{definition} \label{d:hkpsi}
Let $\Psi:[0,\infty) \to [0,\infty)$ be a homeomorphism. For any such $\Psi$, we associate a function $\Phi: (0,\infty) \times (0 ,\infty) \to \bR$ defined by
\[
\Phi(R,t) = \sup_{s > 0}  \left( \frac{R}{s} - \frac{t}{\Psi(s)} \right).
\]	
 We say the MMD space $(\sX,d,\mu,\sE,\sF)$ satisfies heat kernel estimate $\on{HK(\Psi)}$, if
 the heat kernel $p(t,x,y)$ exists and there exists constants $C_1,C_2,C_3,C_4 \in (0,\infty)$ such that
 \begin{align*}
 p(t,x,y) &\ge \frac{1}{\mu \left(B\left(x,\Psi^{-1}(C_1t)\right)\right)} \exp \left(- \Phi(C_2 d(x,y), t)\right), \\ p(t,x,y) &\le \frac{1}{\mu \left(B\left(x,\Psi^{-1}(C_3t)\right)\right)} \exp \left(- \Phi(C_4 d(x,y), t)\right),
 \end{align*} 
 for all $(t,x,y) \in (0,\infty) \times \sX \times \sX$.
\end{definition}
Note that if $\Psi(r)= r^2$, $\on{HK(\Psi)}$ corresponds to Gaussian estimates, and if $\Psi(r)=r^{d_w}$ the estimates are analogous to \eqref{e:usg}, \eqref{e:lsg}.

\subsection{Cable system and Cable process} \label{ss:cable}
A good reference on cable systems and associated Markov processes is \cite{Fol}.
These processes considered in this section can also be viewed as a special (one dimensional) case of diffusions on Riemannian complexes considered in \cite{PS}.
We shall see that a cable system can be viewed as a MMD space and admits a notion of modulus that is compatible with the MMD space.

Consider a connected, locally finite, simple graph  $\bG=(V_\bG,E_\bG)$
 endowed with a \emph{length function} $\ell:E_\bG \to (0,\infty)$. 
 Recall that a simple graph is an undirected graph in which both multiple edges and loops are disallowed.
We view the edges $E_\bG$ as a subset of the two-element subsets of $V_\bG$, \emph{i.e.} $E_\bG \subset \set{ J \subset V_\bG : \abs{J}=2}$.
We define an arbitrary orientation by providing each edge $e \in E_\bG$ with a source $\wh{s}:E_\bG \to V$ and a target $\wh{t}:E_\bG \to V_\bG$ such that $e= \set{\wh{s}(e),\wh{t}(e)}$. We say two vertices  $u,v \in V_\bG$ are \emph{neighbours} if $\set{u,v} \in E_\bG$. We say two distinct edges $e,e' \in E_\bG$ are \emph{incident} if $e \cap e' \neq \emptyset$.

 The {\em cable system}  $\sX=\sX(\bG)$ corresponding to the graph $\bG$ is the topological space obtained
by replacing each edge $e \in E_\bG$ by a copy of the unit interval $[0,1]$, glued together in the obvious way, with the endpoints corresponding to the vertices.
More formally, we define $\sX$ as the quotient space $\left(E_\bG \times [0,1]\right)/\sim$, where $\sim$ is the smallest equivalence relation such that $\wh{t}(e)=\wh{s}(e')$ implies $(e,1) \sim (e',0)$,  $\wh{s}(e)=\wh{s}(e')$ implies $(e,0) \sim (e',0)$, and  $\wh{t}(e)=\wh{t}(e')$ implies $(e,1) \sim (e',1)$.
Here $E_\bG \times [0,1]$ is equipped with the product topology with $E_\bG$ being a discrete topological space.
It is easy to check that the topological space above does not depend on the choice of the edge orientations given by $s,t:E_\bG \to V$. 
This defines the cable system $\sX$ as a topological space equipped with the canonical quotient map
$q: E_\bG \times [0,1] \to \sX$. 
It is easy to check that $\sX$ is  locally compact and separable.

 We denote $\sX_e=q\left( \set{e} \times [0,1] \right)$, so that $\sX= \cup_{e \in E_\bG} \sX_e$. We sometimes abuse notation and abbreviate $q(e,s) \in \sX_e$ by $(e,s)$. There is a canonical injection $i: V_\bG \to \sX$ such that $\wh{s}(e)=v$ implies $i(v)=q(e,0)$ and $\wh{t}(e')=v$ implies $i(v)=q(e',1)$. We abuse notation and abbreviate $i(v)$ by $v$ and therefore view $V_\bG$ as a subset of $\sX$.

We now define a \emph{metric} $d_\ell:\sX \times \sX \to [0,\infty)$ induced by a \emph{length function} $\ell:E_\bG \to (0,\infty)$. 
First, we define a metric $d_{\ell,e}:\sX_e \times \sX_e \to [0,\infty)$ as $d_{\ell,e}((e,s),(e,t))= \ell(e) \abs{s-t}$. By \cite[Corollary 3.1.24 and Exercise 3.2.14]{BBI}, there is a unique maximal metric $d_\ell:\sX \times \sX \to [0,\infty)$ such that $d_\ell(x,y) \le d_{\ell,e}(x,y)$ for all $e \in E_\bG, x,y \in \sX_e$.

We describe some properties of the metric space $(\sX,d_\ell)$. 
 The metric space $(\sX,d_{\ell})$ is a length space such that the metric topology coincides with the quotient topology defined above -- see \cite[Exercise 3.2.14]{BBI}.
We can recover the length function $\ell:E_\bG \to (0,\infty)$ from the metric $d_\ell$ using $\ell(e)= L(\gamma_e)$ where the curve $\gamma_e:[0,1] \to \sX$ given by $\gamma_e(s)= q(e,s)$ \cite[Exercise 3.2.16]{BBI}. 
If $\ell\equiv 1$, then the metric $d_\ell$ restricted to $V_\bG \times V_\bG$ coincides with the  graph distance metric $d_\bG$. 
 This metric space $(\sX,d_\ell)$ is also  called a \emph{metric graph} or \emph{one-dimensional polyhedral complex} -- \cite[Section 3.2.2]{BBI}. 
\begin{remark}
A quicker definition of the metric space $(\sX,d_\ell)$ would be to replace each edge $e$ by an isometric copy of $[0,\ell(e)]$ and glue them in an obvious way at the vertices and consider the induced metric.	
However, we did not follow that approach because it will be important for us to view   $d_\ell$ as a family of metrics on the \emph{same} topological space $\sX$ -- see Lemma \ref{l:modinvariant}. In particular, our definition provides canonical homeomorphism between $(\sX,d_\ell)$ and $(\sX,d_{\wt{\ell}})$ via the identity map, where $\ell,\wt\ell$ are two different length functions.
\end{remark} 

We say $x \in \sX$ is a \emph{vertex} in $\sX$ if $x = q(e,0)$ or $x= q (e,1)$ for some $e \in E_\bG$. We denote the set of vertices in $\sX$ by $\sX_V$.
There is an obvious bijection between vertices of the cable system $\sX_V$ and vertices of the graph $V_\bG$.
For a vertex $x \in \sX_V$, we define the \emph{separation radius} as
\[
r_x= \inf_{y \in \sX_V \setminus \set{x}} d_\ell(x,y) = \min_{\sX_e \ni x} \ell(e).
\]
For non-vertices $x \in \sX \setminus \sX_V$, we define the \emph{separation radius} as
\[
r_x= \min_{ y \in \sX_V: \set{x, y} \subset \sX_e} r_y.
\]
The above minimum is clearly over two vertices in $\sX_V$. Using the local finiteness of the graph, it is easy to see that  $r_x \in  (0,\infty)$ for all $x \in \sX$.

Now, we define a Borel measure $\mu_\ell$ induced by the length function $\ell$.
 The canonical measure $\mu_\ell$ on $\sX$ is defined by $\mu_\ell(\sX_V)=0$ and
 \be \label{e:defmeas}
\mu_\ell\left( \set{e} \times (s,t)\right)=  \ell(e)^2 \abs{s-t}
 \ee
for all $e \in E$ and $0 \le s < t \le 1$.
 The reason for the factor $\ell(e)^2$ above will become apparent in \eqref{e:energy} and Lemma \ref{l:modinvariant}.
 Observe that  $\mu_\ell(\sX_e)=\ell(e)^2$ for all $e\in E_\bG$. 
Further, if $\ell, \wt{\ell}$ are two length functions, then the corresponding measures $\mu_\ell$ and $\mu_{\wt{\ell}}$ are mutually absolutely continuous. 
Therefore one can unambiguously talk about the notion of \emph{almost everywhere} in $\sX$ without even specifying  $\mu_\ell$ or $\ell$. Evidently, $\mu_\ell$ has full support.

 If $\ell\equiv 1$, then the space $(\sX,d_\ell)$ is complete, but completeness need not hold for an arbitrary length function.  
 In general, we denote the completion of $\sX$ with respect to $d_\ell$ by the space $(\ol{\sX},d_\ell)$. 
 We abuse notation and write $d_\ell:\ol{\sX} \times \ol{\sX} \to [0,\infty)$ as the natural extension of the metric on $\sX$, and $\mu_\ell$ to be the Borel measure on $\ol{\sX}$ that extends the measure on $\sX$  such that $\mu_\ell(\ol{\sX} \setminus \sX)=0$. We say $\ell$ is \emph{complete} if $(\sX,d_\ell)$ is complete, \emph{i.e.}, $\sX = \ol\sX$.

\noi \textbf{Notation}. We use $B_\ell(x,r)$ to denote an open ball in $(\sX,d_\ell)$ and $V_\ell(x,r)$ to denote its volume $\mu_\ell\left(B_\ell(x,r)\right)$. If $\ell \equiv 1$, we write $d_\ell,\mu_\ell,B_\ell$ and $V_\ell$ as $d_1,\mu_1,B_1$ and $V_1$ respectively. 

We introduce a notion of \emph{length of the gradient} on the cable system $(\ol\sX,d_\ell,\mu_\ell)$ of a suitable function $f:\ol\sX \to \bR$. We denote the arc-length parametrization of the curve $\gamma_e$ defined above as $\wt{\gamma}_e:[0,\ell(e)] \to \sX$ with $\wt{\gamma}_e(s) = q(e,s/\ell(e))$. We say a function $f:\ol\sX \to \bR$ is \emph{absolutely continuous} if $f$ is continuous and if $f \circ \wt{\gamma}_e: [0,\ell(e)] \to \bR$ is absolutely continuous for each $e \in E_\bG$.   
Clearly, absolute continuity of $f$ implies that $f \circ \wt{\gamma}_e$ is differentiable almost everywhere in $(\ol\sX,\mu_\ell)$.
If $f$ is absolutely continuous, then there is a function $\abs{\nabla_\ell f}:\ol\sX \to [0,\infty)$ such that for each $e \in E_\bG$, we have 
\be \abs{\nabla_\ell f}(\wt{\gamma}_e(s))=\abs{ (f \circ \wt{\gamma}_e)'(s)} \label{e:defgrad} \ee for almost every $s \in [0,\ell(e)]$. Notice that $\abs{\nabla_\ell f}$ is well-defined up to sets of measure zero and 
does not depend on the choice of orientation $\wh{s},\wh{t}:E_\bG \to V$. 

We now relate $\abs{\nabla_\ell f}$ to  upper gradient--see Definition \ref{d:ug}.
By the fundamental theorem of calculus and triangle inequality, $\abs{\nabla_\ell f}$ is an upper gradient of $f$ in $(\ol\sX,d_\ell)$. Further by \cite[Proposition 6.3.3]{HKST}, any upper gradient $\rho$ of $f$ satisfies $\rho \ge \abs{\nabla_\ell f}$ almost everywhere. Therefore  $\abs{\nabla_\ell f}$ is the minimal upper gradient and is unique in the almost everywhere sense.

The \emph{Dirichlet energy} of an absolutely continuous function $f$  on $(\ol\sX,d_\ell)$ is defined by
\be \label{e:energy}
\sD^\ell(f,f)= \int_{\sX} \left(\abs{\nabla_\ell f}(x) \right)^2\, d \mu_{\ell}(x).
\ee
It is easy to check that $\sD^\ell(f,f)$ does not depend on the length function $\ell$. 
We can unambiguously abbreviate $\sD^\ell(f,f)$ as $\sD(f,f)$.

Let $\ell:E_{\bG} \to (0,\infty)$ be a complete length function in a graph $(V_\bG,E_\bG)$.  Let $\Mod^\ell(\cdot)$ denote the modulus of the metric measure space $(\sX,d_\ell,\mu_\ell)$ corresponding to the cable system. 
The modulus $\Mod^\ell$ does not depend on $\ell$ as shown below. The following lemma can be viewed as a conformal invariance of modulus.
\begin{lem}(See \cite[Theorem 7.10]{Hei}) \label{l:modinvariant}
	Let $\bG=(V_\bG,E_\bG)$ be a graph with two complete length functions $\ell,\wt{\ell}:E_\bG \to (0,\infty)$.
		Let $\Gamma$ be a family of curves in  the cable system $\sX$. 
  Then $\Mod^\ell(\Gamma)= \Mod^{\wt{\ell}}(\Gamma)$.
\end{lem}
\proof
Let $\varrho$ be an admissible metric for $\Gam$ in $(\sX,d_\ell,\mu_\ell)$. Then 
\[\wt{\varrho}(x)= \sum_{e \in E} \varrho(x) \mathbf{1}_{\sX_e}(x) \frac{\wt{\ell}(e)}{\ell(e)}\] is an admissible metric for $\Gamma$ in $(\sX,d_{\wt{\ell}},\mu_{\wt{\ell}})$ with $\int_{\sX} \varrho^2\, d\mu_\ell = \int_{\sX} \wt{\varrho}^2\, d\mu_{\wt{\ell}}$. This shows $\Mod^{\wt{\ell}}(\Gam) \le \Mod^{\ell}(\Gam)$, which implies the desired result by symmetry.
\qed\\
By the above Lemma, we can unambiguously abbreviate $\Mod^\ell$ by $\Mod$ for complete length functions $\ell$. 

Now we define a Dirichlet form $(\sE,\sF^\ell)$ associated to a cable system (completed) $(\ol\sX,d_\ell,\mu_\ell)$ on a graph $(V_\bG,E_\bG)$ with length function $\ell$.
Let $\langle \cdot, \cdot \rangle_\ell$ denote the inner product in $L^2(\ol\sX,\mu_\ell)$.
 Let $\sA$ denote the vector space of all absolutely continuous functions $f:\ol\sX \to \bR$ on $(\ol\sX,d_\ell)$, $\sC_c(\ol\sX)$ denote the space of continuous functions with compact support in $\ol\sX$, and $\sW$ denote
\be \label{e:defW}
\sW = \set{ f \in \sA : \sD(f,f) < \infty},
\ee
where $\sD$ denotes the Dirichlet energy in \eqref{e:energy}.

 We define $\sE(f,f)= \sD(f,f)$ for all $f \in \sF^\ell$, where
 the domain of the Dirichlet form $\sF^\ell$ is defined as the closure of  $\sW \cap \sC_c(\ol\sX)$ with respect to the norm $f \mapsto \left(\sD(f,f)+ \langle f,f \rangle_\ell\right)^{1/2}$.
Using properties of Sobolev spaces \cite[Theorem 8.7 and Proposition 8.1]{Bre}, it is easy to check that $(\sE,\sF^\ell)$ is a regular, strongly local Dirichlet form on the locally compact, metric measure space $(\ol\sX,d_\ell,\mu_\ell)$.

Clearly, $\sF^\ell$ is a subspace of $L^2(\ol\sX,\mu_\ell) \cap \sW$. 
Let $\sF_e^\ell$ denote the corresponding extended Dirichlet space.
If $\ell\equiv 1$, we denote $\sF^\ell,\sF^\ell_e$ by $\sF^1,\sF^1_e$ respectively.
The energy measure corresponding to $(\sE,\sF^\ell)$ is easily verified to be $d\Gamma^\ell(f,f) = \abs{\nabla_\ell f}^2 \, d\mu_\ell$.

By Sobolev embedding \cite[Theorem 8.2]{Bre}, every function $f \in \sF_e^\ell$ has a continuous version that is also absolutely continuous. Henceforth, we shall \emph{always} represent every function in $\sF_e^\ell$ (and therefore $\sF^\ell$) by its absolutely continuous version. 
Further, the Dirichlet form $(\sE,\sF^\ell)$ on $L^2(\ol\sX,\mu_\ell)$ is \emph{irreducible} -- see \cite[p. 48]{FOT} for the notion of irreducibility.
\begin{remark} \label{r:timechange}
{\rm
\begin{enumerate}
\item
If $(\sX,d_\ell)$ is complete, then the MMD space $(\sX,d_\ell,\mu_\ell,\sE,\sF^\ell)$ is a time change of $(\sX,d_1,\mu_1,\sE,\sF^1)$. In particular, $(\sX,d_\ell,\mu_\ell,\sE,\sF^\ell)$ is recurrent if and only if $(\sX,d_1,\mu_1,\sE,\sF^1)$ is recurrent. 
\item The above statements do not necessarily hold if $(\sX,d_\ell)$ is not complete. For example, let $(\sX,d_\ell)$ denote the circle packing embedding with straight lines of a one-ended, bounded degree, planar triangulation whose carrier is the open unit disk $\bU$. Then $(\sX,d_\ell,\mu_\ell,\sE,\sF^\ell)$ is recurrent, while $(\sX,d_1,\mu_1,\sE,\sF^1)$ is transient.
\end{enumerate}
}
\end{remark}

Let $(X_t)_{t \ge 0}$ be the Hunt process associated 
with the MMD space $(\ol\sX,d_\ell,\mu_\ell,\sE,\sF^\ell)$, and write for a Borel set $F \subset \ol\sX$,
\be \label{e:defhit}
T_F = \inf\{ t > 0: X_t \in F\}, \, \tau_F = T_{F^\compl}. 
\ee
That is, $T_F$ and $\tau_F$ denote, respectively, the hitting time and exit time of the set $F$.
\subsection{Embedding planar graphs}
An \emph{embedding with straight lines} of a planar graph $\bG=(V_\bG,E_\bG)$ is a map sending the vertices to points in the plane and edges to straight lines connecting the corresponding vertices such that no two edges cross. We define the \emph{carrier of the embedding}, denoted by $\on{carr}(\bG)$ be the union of closed faces of the embedding.

We identify a vertex $v$ with the image in the embedding.
We write $\abs{u-v}$ for the Euclidean distance between the points $u$ and $v$ in the plane.
Any embedding with straight lines defines a length function $\ell:E_\bG \to (0,\infty)$, where $\ell(e)=\abs{u-v}$ for  $e=\set{u,v}$. In this case, we say that $(\sX,d_\ell,\mu_\ell)$, and $(\sX,d_\ell,\mu_\ell,\sE,\sF^\ell)$ are the cable system, and MMD space corresponding to the embedding.

Circle packing naturally leads to an embedding with straight lines which we describe now. 
By drawing the edges as straight lines joining the centers of the corresponding
circles in the circle packing, we obtain an embedding with straight lines  of $\bG$ in
$\bR^2$ such that no two edges cross.  The carrier of a circle packing is defined to be the carrier of the associated embedding with straight lines. 
He and Schramm \cite{HS} showed that a bounded degree, one-ended planar triangulation can be circle packed so that the carrier is either the open unit disk $\bU$ or the entire plane $\bR^2$ depending on whether the simple random walk on $\bG$ is recurrent or transient, respectively.

We recall the notion of \emph{good embedding} from \cite{ABGN}.
\begin{definition} \label{d:good}
Let $D \in (0,\infty)$ and $\eta \in (0,\pi)$.
We say than an embedding with straight lines of a planar graph $\bG=(V_\bG,E_\bG)$ is $(D,\eta)$-good if 
\begin{itemize}
\item[(a)] \emph{No flat angles}. For any face, all the inner angles are at most $\pi - \eta$. In particular, all faces are convex, there is no outer face, and the number of edges in a face is at most $2 \pi/\eta$.
\item[(b)] \emph{Adjacent edges have comparable lengths}. 
For any two adjacent edges $e_1=\set{u,v}$ and $e_2=\set{u,w}$, we have $\abs{u-v}/\abs{u-w} \in [D^{-1},D]$.
\end{itemize}
\end{definition}
If the carrier of an embedding $\on{carr}(\bG)= \bR^2$, then by Hopf--Rinow--Cohn-Vossen theorem  \cite[Theorem  2.5.28]{BBI} the corresponding cable system $(\sX,d_\ell,\mu_\ell)$ is complete.

  The ring lemma of Rodin and Sullivan shows that circle packing induces a good embedding of bounded degree triangulations \cite[p. 352]{RS}. More generally, the following lemma provides a large family of planar graphs that admit a good embedding.
  
  \begin{lemma}(\cite[p. 352]{RS} and \cite[Corollary 4.2]{HN}) \label{l:good}
  Let $\bG$ be an one-ended, bounded degree, simple, 3-connected, planar graph, such that its planar dual  $\bG^\dagger$ is also a bounded degree graph. Then $\bG$ admits a good embedding in the sense of Definition \ref{d:good}.
  \end{lemma}

\section{Annular quasi-convexity at large scales} \label{s:aqc}
The goal of this section is to obtain a geometric consequence of modulus estimates. Roughly speaking, Poincar\'e inequality implies a lower bound on the modulus and upper bound  on capacity implies a upper bound  on the modulus. Estimates on modulus have useful geometric consequences. One of the geometric properties that will be studied in this section is the notion of linear local connectivity, which in a sense says every annulus $B(x,2r) \setminus B(x,r)$ is `well connected' -- see Definition \ref{d:llc}.

To obtain bounds on the modulus, we show an analogue of the `Dirichlet's principle' for modulus and provide a probabilistic formula to compute modulus.
\begin{lemma} \label{l:dp} 
	Let $(\sX,d_\ell,\mu_\ell,\sE,\sF^\ell)$ denote a complete cable system.
	Let $E,F$ be nonempty, disjoint sets such that $\dist_\ell(E,F)>0$, $\diam_\ell(E)<\infty$, and $\bP_x(T_{E \cup F}  < \infty)=1$ for all $x \in \sX$. Then $\Mod(E,F) = \int_{\sX} d \Gamma^\ell(u,u) = \int_{\sX} \abs{\nabla_\ell u}^2\, d\mu_\ell$, where $u(x)= \bP_x(T_E < T_F) \in \sF_{\loc}^\ell(\sX)$, and $T_E,T_F$ denotes the hitting times as defined in \eqref{e:defhit}.
\end{lemma}
\proof
First, we  show that every ball $B(x,r)$ has finite measure.
By \cite[Proposition 2.5.22]{BBI}, the closure $\ol{B(x,r)}$ is compact.
Let $\mathfrak{U}$ denote the open cover of $\sX$ formed by the interior of edges  and small balls around vertices in $\sX$ given by
\be \label{e:cover}\mathfrak{U}=\set{\sX_e^\circ : e \in E} \cup \set{B_\ell( x,r_x/4): x \in \sX_V}. \ee
By compactness, we see that that  $\ol{B(x,r)}$ intersects only finitely many edges and therefore has finite measure.

Let $\Gamma(E,F)$ denote the family of curves joining $E$ and $F$. 
Let $B$ be a ball that contains $\set{y \in \ol\sX: \dist_\ell(y,E) \le \dist_\ell(E,F)}$.
Consider
\[
f(x) = \frac{1}{\dist_\ell(E,F)}\one_{B}(x).  
\]
Then $f$ is an admissible metric for $\Gamma(E,F)$ which is also in $L^2(\ol\sX,\mu_\ell)$. Therefore $\Mod(E,F)$ is finite and it suffices to restrict our attention to admissible metrics in $L^2(\ol\sX,\mu_\ell)$.

Clearly, the space of admissible metrics for $\Gamma(E,F)$ is convex. Therefore by uniform convexity of $L^2(\ol\sX,\mu_\ell)$, there is at most one \emph{optimal} admissible metric $\rho$ for $\Gamma(E,F)$ such that $\Mod(E,F)= \int \rho^2 \, d\mu_\ell$. This establishes \emph{uniqueness of optimal admissible metric}. Next, we proceed to show \emph{existence} of such a metric.

To prove existence, we will show that the space of admissible metrics of $\Gamma(E,F)$ that are in $L^2(\ol\sX,\mu_\ell)$ forms a closed subset of $L^2(\ol\sX,\mu_\ell)$.
Let $\rho_n$ be a sequence of admissible metrics for $\Gamma(E,F)$ that converges to $\rho$ in $L^2(\ol\sX,\mu_\ell)$. We now show that $\rho$ is also admissible for $\Gam(E,F)$.
Let $\gamma:[a,b] \to \sX$ be a curve joining $E$ and $F$. By removing loops if necessary, it suffices to assume that $\gamma$ is simple. 
Let $A$ denote the union of the edges $\sX_e$ such that $\sX_e \cap \gamma([a,b]) \neq \emptyset$. Since $\gamma([a,b])$ is compact, by using the open cover $\mathfrak{U}$ in \eqref{e:cover}, we obtain that $A$ is a finite union of edges. Therefore, there exists a constant $C_\gamma \in (0,\infty)$ such that
\begin{align*}
\abs{ \int_{\gamma} (\rho_i - \rho)\, ds} &\le \int_{\gamma} \abs{\rho_i-\rho}\,ds \\
 &\le C_\gamma \int_{A}  \abs{\rho_i-\rho}\, d\mu_\ell \le C_\gamma \left( \int_{A}  \abs{\rho_i-\rho}^2\, d\mu_\ell \right)^{1/2} \left(\mu_\ell(A) \right)^{1/2}.
\end{align*}
The above estimate shows that $\int_{\gamma} \rho_i \,ds \to \int_{\gamma} \rho_i \,ds$. Therefore $\rho$ is admissible.

Let $\rho_i \in L^2(\ol\sX,\mu_\ell), i \in \bN$ be a sequence of admissible metrics such that
\be 
\int_{\ol\sX} \rho_i^2 \, d\mu_\ell \le \Mod(E,F) + \frac{1}{i}.
\ee
By Banach-Saks theorem  and by passing through a subsequence if necessary (and denote the subsequence again by $\rho_i$), we can assume that the Ces\`aro means $\wt{\rho}_i = i^{-1} \sum_{j=1}^i \rho_i$ converge in $L^2(\sX,\mu_\ell)$ to $\rho \in L^2(\sX,\mu_\ell)$.
By the triangle inequality, 
 $\Mod(E,F)= \int_{\sX} \rho^2 \, d\mu_\ell$ and $\rho$ is admissible for $\Gamma(E,F)$. This completes the proof of existence of an optimal admissible metric $\rho$.

Let $\Gamma_x$ denote the family of curves that join $x$ to $F$.
Define 
\be \label{e:dp0}
u(x) = \inf_{\gamma \in \Gamma_x} \int_\gamma \rho\,ds.\ee
Clearly $\restr{u}{F} \equiv 0$ and
by admissibility of $\rho$, we have $\restr{u}{E} \ge 1$. 
Further, by Sobolev embedding \cite[Theorem 8.2]{Bre}, the function $u$ is absolutely continuous on $\sX$. 
Any upper gradient of $u$ is admissible for $\Gamma(E,F)$ \cite[Proof of Propositon 2.17]{HK}.
Since $\rho$ is optimal, $\rho$ is the unique minimal upper gradient of $u$ \cite[Theorem 6.3.20]{HKST} and hence $\rho = \abs{\nabla_\ell u}$.

Next, we show that $u \in \sF^\ell_e$.
Similarly,  any upper gradient of  $\wt{u}= (u \vee 0) \wedge 1$  is admissible for $\Gamma(E,F)$. The minimal upper gradient $\abs{\nabla_\ell \wt{u}}$ of $\wt{u}$ satisfies $\abs{\nabla_\ell \wt{u}} \le \abs{\nabla_\ell u}$ almost everywhere \cite[Proposition 6.3.23]{HKST}. By the optimality of $\rho=\abs{\nabla_\ell u}$, admissibility of $\abs{\nabla_\ell \wt{u}}$, $\abs{\nabla_\ell \wt{u}} \le \abs{\nabla_\ell u}$, and  uniqueness of optimal metric, we have $\abs{\nabla_\ell \wt{u}} \equiv \abs{\nabla_\ell u}$. Hence by replacing $u$ by $\wt{u}$ if necessary, henceforth we shall assume $0 \le u \le 1$.

We will show that 
\be  \label{e:dp1}
f \wedge g \in \sF^\ell, \hspace{3mm} \mbox{ for all non-negative functions $f \in \sF^\ell, g \in \sW$,}
\ee
where $\sW$ is as defined in \eqref{e:defW}.
By regularity, there exists a sequence $f_n \in \sF, n \in \bN$, such that
$f_n \to f$ in the norm $h \mapsto \left(\sE(h,h) + \langle h,h\rangle_{\ell} \right)^{1/2}$. By passing through a subsequence if necessary, we can assume that $f_n \to f$ almost everywhere.
Clearly $f_n \wedge g \in \sW \cap \sC_c(\sX) \subset \sF$, and 
\be \label{e:dp2}
\liminf_{n \to \infty} \sE(f_n \wedge g, f_n \wedge g) \le \int_{\sX} \abs{\nabla_\ell g}^2\,d\mu_\ell + \liminf_{n \to \infty} \sE(f_n,f_n) = \int_{\sX} \abs{\nabla_\ell g}^2\,d\mu_\ell +  \sE(f,f) < \infty.
\ee
Using \eqref{e:dp2} and Lemma \ref{l:extd}, we have $f \wedge g \in \sF^\ell_{\loc}$. Since $ \abs{f \wedge g} \le g$, we have $f \wedge g \in L^2(\sX,\mu_\ell)$, which along with Lemma \ref{l:extd} establishes \eqref{e:dp1}. 

Next, we check that $u \in \sF^\ell_{\loc}(\sX)$ using \eqref{e:dp1}. Since $u \in \sW$,  it is continuous and hence bounded on any compact subset $K$. By the regularity of the Dirichlet form, there exists a non-negative function $f \in \sF \cap \sC_c(\sX)$ such that $f> u$ on a neighbourhood of $K$. Therefore $u \wedge f \in \sF$ such that $u \wedge f =u$ on $K$. Since $K$ was arbitrary compact set, we have $u \in \sF_{\loc}^\ell(\sX)$.


%

Similarly, if $v$ is absolutely continuous in $\sX$ with $\restr{v}{E} \ge 1, \restr{v}{F} \equiv 0$, then $\abs{\nabla_\ell v}$ is admissible for $\Gamma$.
Hence for $\phi \in \sC_c( \sX \setminus (E \cup F)) \cap \sF$, the function
\[
t \mapsto \sE(u+t \phi,u+t\phi)
\]
is minimal at $t=0$. Therefore $\sE(u,\phi)=0$ for all  $\phi \in \sC_c( \sX \setminus (E \cup F)) \cap \sF^\ell$.
This implies that $u \in \sF_{\loc}^\ell(\sX)$ is harmonic in $\sX \setminus (E \cup F)$.

Since $u$ is harmonic in $(E \cup F)^c$, by Remark \ref{r:harm} $M_t= u(t \wedge  T_{E \cup F})$ is a bounded martingale. Since $\bP^x(T_{E \cup F} < \infty)=1$,  the martingale $M_t$ converges almost surely to $\one_{\set{T_E < T_F}}$.
By the optional stopping theorem, we have
$u(x)= \bP^x(T_E < T_F)$. \qed
\begin{remark} \label{r:modcap}
	The above proof shows that $\Mod(E,F) \le \Cap_{F^c}(E)$. 
The reverse inequality is not true in general because the function $u$ in \eqref{e:dp0} need not be in the extended Dirichlet space.	
\end{remark}

The following definition is a large scale variant of \cite[3.12]{HK} -- see also \cite{Kor} for the related notion of annular quasi-convexity.
\begin{definition}[Annular quasi-convexity at large scales] \label{d:llc}
We say that a metric space $(\sX,d)$ is annular quasi-convex at large scales if there exists $R_L,C_L>0$ such that
for all $x \in \sX$, for all $r  \ge R_L$, and for all $y,z \in B(x,2r)\setminus B(x,r)$ there exists a curve $\gamma:[a,b] \to B(x,C_L r) \setminus B(x,r/C_L)$ that connects $y=\gamma(a)$ and $z=\gamma(b)$.
\end{definition}
\begin{prop} \label{p:llc}
	Let $\bG=(V_\bG,E_\bG)$ be a graph of polynomial volume growth with volume growth exponent $d \ge 2$. 
	Assume that the corresponding cable system $(\sX,d_1,\mu_1,\sE,\sF^1)$ satisfies the capacity upper bound  \hyperlink{cap}{$(\on{cap}_{\le})$},  and Poincar\'e inequality 
	\hyperlink{pi}{$\on{PI}(\Psi)$}, with $\Psi(r)= r^2 \vee r^d$.
	Then $(\sX,d_1)$ is annular quasi-convex at large scales.
\end{prop}
\proof 
We follow the argument in \cite[Lemma 3.17]{HK}. 
Since $\bG$ has polynomial growth with exponent $d$, the cable system $(\sX,d_1,\mu_1)$ satisfies the following volume estimate: there exists $C>1$ such that for all $x \in \sX, r \ge 1$, 
\be \label{e:llc0}
C^{-1} r^d \le V_1(x,r) \le C  r^d.
\ee

Let $y,z \in B_1(x,2r) \setminus B_1(x,r)$ and define sets $E=\overline {B_1(y,r/8)}$ and $F=\overline{B_1(z,r/8)}$. 
If $d(y,z)\le r/4$,  then the geodesic from $y$ to $z$ does not intersect $B_1(x,r/2) \cup B_1(x, 3r)^c$. Therefore, we may assume $d_1(y,z) > r/4$. Hence $E,F$ are disjoint, non-empty closed sets with $\dist_1(E,F) > 0$.
By Lemma \ref{l:series}, one easily checks that the MMD space is recurrent.
By Lemma \ref{l:dp}, there is a function $u \in \sF_{\loc}^1(\sX)$ with $\Mod(E,F) = \int_{\sX} \abs{\nabla_1 u}^2 \, d\mu_{1}$, $\restr{u}{E}\equiv 1, \restr{u}{F} \equiv 0$.

Let $M \ge 1$ be the constant in \hyperlink{pi}{$\on{PI}(\Psi)$}.
By using \hyperlink{pi}{$\on{PI}(\Psi)$} for $u$ in a ball $B_1(x, 3Mr)$, and using \eqref{e:llc0}, we obtain the following lower bound on $\Mod(E,F)$: there exists $c>0$ such that for all $x \in \sX, r \ge 1$,  and $y,z \in B_1(x,2r) \setminus B_1(x,r)$ with $d(y,z)  > r/4$, we have
\be \label{e:llc1}
\Mod(E,F) \ge c.
\ee

Let $C>4$ be a constant, whose value will be momentarily determined.
Let $\Gamma_1$  be the family of curves joining $E \cup F$ and $B_1(x,r/C)$, $\Gamma_2$ be the family of curves joining $E \cup F$ and $\sX \setminus B_1(x,Cr)$, and $\Gamma_3$ be the family of curves joining $E$ and $F$ in $B_1(x,Cr) \setminus B_1(x,r/C)$.
By the basic properties of modulus \eqref{e:bpm2}-\eqref{e:bpm4}, we have
\begin{align} \label{e:llc2}
\lefteqn{\Mod(E,F) \le \sum_{i=1}^3 \Mod(\Gamma_i) }\nonumber \\
&\le \Mod(\Gamma_3) + \Mod(B_1(x,r/C),\sX \setminus B_1(x, r/2)) + \Mod(B_1(x,3r)), \sX \setminus B_1(x,Cr)).
\end{align}
For any function $u$ that is admissible in the definition of $\Cap_{B_1(x,R)}(B_1(x,r))$, its gradient $\rho=\abs{\nabla_\ell u}$ is an admissible metric for $\Mod(B_1(x,r),\sX \setminus B_1(x,R)$.
Therefore,for all $x \in \sX$, $0<r < R$,
\be  \label{e:llc2a}
 \Mod(B_1(x,r),\sX \setminus B_1(x,R)) \le \Cap_{B_1(x,R)}(B_1(x,r))
 \ee
By Lemma \ref{l:series} and \eqref{e:llc2a}, there exists $C,R > 1$ such that for all $r \ge R$, $x \in \sX$, we have
\be \label{e:llc3}
\Mod(B_1(x,r/C),\sX \setminus B_1(x, r/2)) + \Mod(B_1(x,3r), \sX \setminus B_1(x,Cr)) \le c/2.
\ee
Combining \eqref{e:llc1}, \eqref{e:llc2}, \eqref{e:llc3}, we have $\Mod(\Gamma_3) \ge c/2 >0$, for all $r \ge R$, $x \in \sX$, and so by \eqref{e:bpm1}, $\Gamma_3 \neq \emptyset$.
\qed

\section{Loewner property of good embeddings} \label{s:loew}
The goal of this section is prove that a good embedding with carrier $\bR^2$ satisfies the Loewner property.
We recall the notion of Loewner space \cite[Chapter 8]{Hei}.

For a metric space, $(\sX,d)$, we  denote distance between two sets or distance between a point and a set by $\dist(E,F) = \inf_{x \in E, y\in F} d(x,y)$ and $\dist(x,E)=  \dist(\set{x},E)=\inf_{y \in E} d(x,y)$ respectively. We denote by $\Delta(E,F)$, the \emph{relative distance} between $E$ and $F$ as 
\[
\Delta(E,F) := \frac{\dist(E,F)}{\diam(E) \wedge \diam(F)}.
\]
For cable systems, the notations $\dist_\ell, \Delta_\ell,\dist_1,\Delta_1$ are self-explanatory.
\begin{definition} \label{d:loew}
	{\rm 
By a \emph{continuum}, we mean a connected, compact set consisting of more than one point.	
We call $(\sX,d,\mu)$ a \emph{Loewner space}, if there is a function $\phi:(0,\infty) \to (0,\infty)$ such that
\be  \label{e:loew}
\Mod(E,F) \ge \phi(t),
\ee
whenever $E$ and $F$ are disjoint continua and $t$ satisfies
\[
t \ge \Delta(E,F).
\]
We say that a space satisfies the Loewner property if it is a Loewner space.
}
\end{definition}
An important aspect of the above definition is that the quantity $\Delta(E,F)$ is `scale-invariant', \emph{i.e.}, $\Delta(E,F)$ does not change if we replace the metric $d$ by $\lambda d$ for some $\lambda > 0$. 
Therefore, the Loewner property can be interpreted as a scale invariant lower bound on the modulus.
The Loewner property is introduced in \cite{HK} and is motivated by Loewner's work on such lower bounds on modulus in the Euclidean space \cite{Loe}.


The main results of \cite{ABGN} show that good embeddings inherit several properties from $\bR^2$; for example Gaussian heat kernel bounds, the volume doubling property, and the Poincar\'e inequality.
The main tool for showing the Loewner property is the Poincar\'e inequality.
The following definition will be used in the proof of the Poincar\'e inequality for good embeddings.
\begin{definition}[Remote balls]
	{\rm
	Let $(\ol{\sX},d)$ be a metric space and $\Omega \subsetneq \ol\sX$ be open. Let $\epsilon > 0$.  We say that a ball $B(x,r)$ is \emph{$\epsilon$-remote} in $\Omega$, if $B(x,r) \subset \Omega$ and $r \le \epsilon d(B(x,r), \sX \setminus \Omega)$. 
	If the value of $\epsilon$ is unimportant, we drop the parameter $\epsilon$ and say $B(x,r)$ is a remote ball.}
\end{definition}

We state a slight generalization of some results in \cite{ABGN}. 
\begin{theorem}(See \cite[Lemma 3.3 and Theorem 3.4]{ABGN}) \label{t:abgn}
	Let $(\overline{\sX},d_\ell,\mu_\ell,\sE,\sF^\ell)$ denote the cable system corresponding to a good embedding of a planar graph with either $\on{carr}(\bG)= \bU$ or $\on{carr}(\bG)=\bR^2$. Then   $(\overline{\sX},d_\ell,\mu_\ell,\sE,\sF^\ell)$ satisfies 
	the volume growth estimate 
\be \label{e:volgood}
\exists C>1: \hspace{3mm} C^{-1} r (r \vee r_x) \le V_\ell(x,r) \le C r (r \vee r_x) \hspace{2mm} \mbox{ for all $x \in \sX, 0<r<\diam_\ell(\sX)/2$},
\ee	
	 and the Poincar\'e inequality  \hyperlink{pi}{$\on{PI}(2)$}.
\end{theorem}

\proof 
The estimate  \eqref{e:volgood} is essentially contained in \cite{ABGN}. In the case, $\on{carr}(\bG)=\bR^2$, one easily checks that the restriction $r <1$ in the statement of \cite[Lemma 3.3]{ABGN} is unnecessary.
The Poincar\'e inequality \hyperlink{pi}{$\on{PI}(2)$} for the case $\on{carr}(\bG)=\bR^2$ is contained in \cite[Theorem 3.4]{ABGN}. 

Although we do not need \hyperlink{pi}{$\on{PI}(2)$} for the case  $\on{carr}(\bG)=\bU$, we provide a proof below.
For the case $\on{carr}(\bG)=\bU$, the Poincar\'e inequality in \cite[Theorem 3.4]{ABGN} is proved only for remote balls in $\sX \subsetneq \ol\sX$.
However, as mentioned in the beginning of \cite[Section 3.2]{ABGN}, the weak Poincar\'{e} inequality implies the strong Poincar\'{e} inequality for all remote balls in $\sX \subset \overline{\sX}$. This is due to Jerison using a Whitney covering argument \cite{Jer86} (see also \cite[Theorem 4.18]{Hei} and \cite[Corollary 5.3.5]{Sal02}).

For inner uniform domains satisfying the volume doubling property, Poincar\'{e} inequality for remote balls implies Poincar\'{e} inequality for all balls \cite[Theorem 3.13]{GyS}. This was established using a Whitney covering in \cite{GyS}. Hence \hyperlink{pi}{$\on{PI}(2)$} for the case $\on{carr}(\bG)=\bU$  follows from \cite[Theorem 3.4 and Lemma 2.6]{ABGN}, \cite[Theorem 3.13]{GyS} and \eqref{e:volgood}.
\qed
\begin{remark}
	{\rm 
For the case $\on{carr}(\bG)= \bU$, the \hyperlink{pi}{$\on{PI}(2)$} stated above implies extensions of \cite[Theorem 1.5 and Theorem 3.6]{ABGN}, where the assumption of remote balls can be relaxed to all balls that are proper subsets of $\ol\sX$. Such an extension was already conjectured in \cite[end of p.1961]{ABGN}. Indeed, the authors of \cite{ABGN} even propose to prove such a generalization in the future by a suitable modification of the  graph.
		Our approach is different from their proposed one and does not require modifying the graph.
	}
\end{remark}

  Recall that the \emph{Hausdorff $s$-content} of a set $E$ in a metric space $(\sX,d)$ is the number
  \[
  \sH^\infty_s(E) = \inf \sum_{i} r_i^s,
  \]
  where the infimum is taken over all countable covers of the set $E$ by balls $B_i$ of radius $r_i$. If $E$ is a continuum in a length space $\sX$, then the Hausdorff $1$-content is comparable to its diameter as
  \be \label{e:diam1content}
  \frac{1}{2} \diam(E) \le \sH^\infty_1(E) \le \diam(E).
  \ee
The upper bound on $\sH^\infty_1(E)$ is easily obtained by covering $E$ using a single ball while the lower bound is  contained in  \cite[proof of Lemma 2.6.1]{BBI}.

\begin{theorem} \label{t:loewner}
Let $(\sX,d_\ell,\mu_\ell)$ be the  cable system corresponding to a good embedding of a planar graph with carrier $\bR^2$. Then  $(\sX,d_\ell,\mu_\ell)$ satisfies the Loewner property. 	
\end{theorem} 
\proof
We follow the argument in \cite[Theorem 5.9]{HK} at  large scales, and then make some essential modifications using local regularity to handle smaller scales.

 
Note that the cable process $(X_t)$ associated to the MMD space $(\sX,d_\ell,\mu_\ell,\sX,\sF^\ell)$ satisfies Gaussian heat kernel upper and lower bounds by Theorem \ref{t:abgn} and \cite[Corollary 4.2 and  4.10]{St}. Integrating this heat kernel estimate in time, we obtain that every Borel set of positive measure is hit infinitely often almost surely. Now using \cite[Theorem 2.1]{Fol}, we conclude that every singleton set is hit almost surely by the cable process from any starting point, that is
\be \label{e:lw1}
\bP^x(T_{\set{y}} < \infty) =1,
\ee
for all $x , y \in \sX$, where 
$T_{\set{y}}$ is as defined in \eqref{e:defhit} for the cable process associated to the MMD space $(\sX,d_\ell,\mu_\ell,\sX,\sF^\ell)$.

Let $t >0$, and
let $E,F \subset  \sX$ be disjoint continua such that $\dist_\ell(E,F) \le t \diam_\ell(E) \wedge \diam_\ell(F)$.

If $\dist_1(E,F)\le 2$, consider a curve $\gamma$ joining $E$ and $F$ in $(\sX,d_1,\mu_1)$ of length at most $2$. Now, by applying the classical Poincar\'e inequality to the function $u$ defined in Lemma \ref{l:dp} on a bounded interval in $\bR$ (corresponding to $\gamma$ in $(\sX,d_1,\mu_1)$), and the conformal invariance of modulus in Lemma \ref{l:modinvariant},  we obtain the following lower bound on modulus: \[\Mod(E,F) \ge \frac{1}{\dist_1(E,F)}, \hspace{4mm} \mbox{ for all pairs of disjoint continua  $E,F \subset \sX$}. \]
Therefore, without loss of generality, we will assume that $\dist_1(E,F) \ge 2$. 

By Lemma \ref{l:dp} and \eqref{e:lw1}, there exists a continuous function $u \in \sF^\ell_{\loc}(\sX)$ such that 
$\restr{u}{E} \equiv 1,\restr{u}{F} \equiv 0$, $u(x)=\bP^x(T_E < T_F)$ and 
\be \label{e:lw1a}
\Mod(E,F)= \int_{\sX} \abs{\nabla_\ell u}^2 \, d \mu_\ell.
\ee
 We need the following gradient estimate: there exists $C_2>0$ such that 
\be \label{e:lw2}
\abs{\nabla_\ell u}(x) \le C_2/r_x
\ee
for almost every $x \in \sX$, where $r_x$ denotes the separation radius. Since $\sX_V$ has measure zero and $\abs{\nabla_\ell u }(x) = 0$  for almost every $x \in E\cup F$, it suffices to consider  $x \in \sX \setminus (\sX_V \cup E \cup F)$.

Every  $x \in \sX \setminus (\sX_V \cup E \cup F)$  belongs to an unique $\sX_e$ for some edge $e$.
We consider two cases depending on whether or not $\sX_e \cap (E \cup F)$ is empty.
If  $\sX_e \cap (E \cup F) = \emptyset$, since the value of $u$ at endpoints of $\sX_e$ differ by at most $1$ and $u$ is linear in the edge $\sX_e$, we have $\restr{\abs{\nabla_\ell u}}{ \sX_e} \le 1/{\ell(e)}$, which immediately implies \eqref{e:lw2}, where $C_2$  depends only on the constants associated with good embedding.

If $x \in \sX_e \cap (\sX \setminus (\sX_V \cup E \cup F))$ is such that  $\sX_e \cap (E \cup F) \neq \emptyset$,  then using $\dist_1(E,F) \ge 2$ we have that $\sX_e$ intersects exactly one of the sets $E$ or $F$.
By symmetry, it suffices to consider the case $\sX_e \cap E \neq \emptyset$. Consider the vertex $v \in \sX_V \cap \sX_e$ such that  $v$ and $x$ belong to the same connected component $I_{x}$ of $\sX_e \setminus E$. Consider the cable process starting at the vertex $v$, exiting the star shaped set $I_x \cup \left(  \cup_{e: \sX_e \ni v} \sX_e\right)$.
By the harmonic measure of this star shaped set from \cite[Theorem 2.1]{Fol} and using $u(y)=\bP^y(T_E < T_F)$ from Lemma \ref{l:dp}, we obtain the gradient estimate \eqref{e:lw2} in this case as well.

We use the notation 
\be \label{e:lw3}
u_{x,r} = \frac{1}{V_\ell(x,r)}\int_{B_\ell(x,r)} u \,d\mu_{\ell}.
\ee
By the gradient estimate \eqref{e:lw2} and the fundamental theorem of calculus, there exists $C_3 >1$ such that for all $x \in \sX, 0<r \le C_3^{-1} r_x$,
we have
\be \label{e:lw4}
\abs{u(x)-u_{x,r}} \le  \frac{1}{V_\ell(x,r)}\int_{B_{\ell}(x,r)} \abs{u(x) - u(y)} \, \mu_\ell(dy) \le \frac{1}{10}.
\ee 

Without loss of generality, assume
\bes
\diam_\ell(E) \le \diam_\ell(F).
\ees
Since $F$ is compact, we can choose $y_0 \in F$ such that $\dist_\ell(E,F)=\dist_\ell(E,\set{y_0})$.
Let $d_F$ denote the geodesic metric induced by $d_\ell$ on $F$ and let $\ol{B_F}(y,r), y \in F,r >0$ denote the corresponding closed balls. By the continuity of the function $s \mapsto \diam_\ell(F \cap \ol{B_F}(y_0,s))$, there exists $s_0>0$ such that $\diam_\ell(F \cap \ol{B_F}(y_0,s_0))=\diam(E)$. 
Replacing $F$ by $F \cap \ol{B_F}(y_0,s_0)$ if necessary and using \eqref{e:bpm2}, we assume that $E,F$ are disjoint continua such that $\dist_\ell(E,F)\le t\diam_\ell(E)=t\diam_\ell(F)$ and $\dist_1(E,F) \ge 2$.
 Since the embdedding is good, the above conditions imply the following bound on the separation radius:  there exists $C_4 >0$ (depending only on $t$ and the constants $D,\eta>0$ associated to the good embedding in Definition \ref{d:good}) such that
\be \label{e:lw5}
r_x \le C_4 \diam_\ell(E) = C_4 \diam_\ell(F) \hspace{2mm} \mbox{ for all $x \in E \cup F$.}
\ee

Fix $R= (3+t)\diam(E)$. By triangle inequality we have $E \cup F \subset B_\ell(x,R)$ for any $x \in E \cup F$.

The proof splits into two cases, depending on whether or not there are points $x \in E$ and $y \in F$ so that neither 
\[
\abs{u(x) - u_{x,R}} \hspace{3mm} \mbox{nor} \hspace{3mm} \abs{u(y) -u_{y,2R}}
\] 
exceeds $\frac{1}{5}$. If such points $x \in E, y \in F$ can be found, then
\[
1 \le \abs{u(x)-u(y)} \le \frac{1}{5} + \abs{u_{x,R} - u_{y,2R}} + \frac{1}{5}.
\]
Therefore, we have
\begin{align*}
\frac{3}{5}\le \abs{u_{x,R} - u_{y,2R}}   &\le \frac{C}{R^2} \int_{B_\ell(y,2R)} \abs{u - u_{y,5R}}\, d\mu_\ell \\
&\le  \frac{C}{R} \left( \int_{B_\ell(y,2R)} \abs{u - u_{y,2R}}^2\, d\mu_\ell \right)^{1/2}\\
& \le \frac{C }{R} \left( R^2 \int_{B_\ell(y,2KR)} \abs{\nabla_\ell u}^2\, d\mu_\ell \right)^{1/2}   \le C \left( \sE(u,u) \right)^{1/2},
\end{align*}
which along with \eqref{e:lw1a} implies Loewner property at large scales. In the above display, we used $B_\ell(x,R) \subset B_\ell(y,5R)$, \eqref{e:lw5} and \eqref{e:volgood} in the first line, 
Cauchy-Schwarz inequality, \eqref{e:lw5} and \eqref{e:volgood} in the second line, Poincar\'e inequality  \hyperlink{pi}{$\on{PI}(2)$} in the third line, and \eqref{e:lw5}, \eqref{e:volgood} in the final line.

The second alternative, by symmetry, is
\be \label{e:lw6}
\abs{u(x) - u_{x,R} } \ge \frac{1}{5} \hspace{3mm} \mbox{ for all $x \in E$.}
\ee
For each $x \in E$, let $i_x \in \bN$ be the unique integer such that
\be \label{e:lw7}
(2 C_3)^{-1} r_x < 2^{-i_x} R \le  C_3^{-1} r_x, 
\ee
so that by \eqref{e:lw4}, and \eqref{e:lw6},
we have
\be \label{e:lw8}
\abs{u_{x,2^{-i_x}R} - u_{x,R}} \ge \frac{1}{10} \hspace{3mm} \mbox{ for all $x \in E$}.
\ee
Using Cauchy-Schwarz inequality, Poincar\'e inequality \hyperlink{pi}{$\on{PI}(2)$}, \eqref{e:lw7} and \eqref{e:volgood}, we obtain the following estimate: for all $x \in E$,
\begin{align*}
1 &\le C  \sum_{j=0}^{i_x-1} \abs{u_{x,2^{-j}R} - u_{x,2^{-j-1}R}}  \le C  \sum_{j=0}^{i_x-1} \frac{1}{V_\ell(x,2^{-j}R)}\int_{B_\ell(x,2^{-j}R)} \abs{u - u_{x,2^{-j}R}} \, d\mu_\ell\\
& \le C  \sum_{j=0}^{i_x-1} \left( \frac{1}{V_\ell(x,2^{-j}R)}\int_{B_\ell(x,2^{-j}R)} \abs{u - u_{x,2^{-j}R}}^2 \, d\mu_\ell \right)^{1/2}\\
&\le  C  \sum_{j=0}^{i_x-1} \left( \frac{(2^{-j}R)^2}{V_\ell(x,2^{-j}R)}\int_{B_\ell(x,2^{-j}MR)} \abs{ \nabla_\ell u}^2 \, d\mu_\ell \right)^{1/2} \\
&\le C  \sum_{j=0}^{i_x-1} \left(\int_{B_\ell(x,2^{-j}MR)} \abs{ \nabla_\ell u}^2 \, d\mu_\ell \right)^{1/2}.
\end{align*}
Therefore, if
\[
\int_{B_\ell(x,2^{-j}MR)} \abs{ \nabla_\ell u}^2 \, d\mu_\ell \le \epsilon 2^{-j},
\]
for some $\epsilon >0$ and for every $x \in \sX, 0 \le j \le i_x -1$, we have that
\[
1 \le C \epsilon^{1/2} \sum_{j=0}^{i_x-1} 2^{-j} \le C \epsilon^{1/2}.
\]
Therefore, for each $x \in E$, there exists an integer $j_x$ with $0 \le j_x \le i_x-1$, such that
\be \label{e:lw9}
\int_{B_\ell(x,2^{-j_x}MR)} \abs{ \nabla_\ell u}^2 \, d\mu_\ell \ge \epsilon_0 2^{-j_x},
\ee
for some small enough $\epsilon_0$ depending only on the constants associated with the definition of good embedding. By the $5B$-covering lemma (see \cite[Theorem 1.2]{Hei} or \cite[p. 60]{HKST}), and the separability of $\sX$, there exists a  countable  family of pairwise disjoint balls $B_k=B_\ell(x_k,2^{-j_{x_k}}MR)$, such that
\be \label{e:lw10}
E \subset \bigcup_{k} B_\ell(x_k,2^{-j_{x_k}}5 M R),
\ee
and, by \eqref{e:lw9}, such that
\be \label{e:lw11}
\diam\left(B_\ell(x_k,2^{-j_{x_k}}5 M R)\right)\le 2^{-j_x+5}MR  \le C R \int_{B_\ell(x_k,2^{-j_{x_k}} M R)}  \abs{ \nabla_\ell u}^2 \, d\mu_\ell.
\ee
Hence, by \eqref{e:diam1content}, \eqref{e:lw10}, \eqref{e:lw11} and the fact that $B_k$'s are disjoint
\begin{align*}
 R/8= \frac{1}{2} \diam(E) &\le \sH_1^\infty(E) \le \sum_k \diam\left(B_\ell(x_k,2^{-j_{x_k}}5 R)\right)\\
 & \le  C R \sum_{k} \int_{B_k}  \abs{ \nabla_\ell u}^2 \, d\mu_\ell \le CR \Mod(E,F).
\end{align*}
This completes the proof of Loewner property.
\qed
\section{Quasisymmetry of good embeddings}  \label{s:qs}

\begin{theorem} \label{t:qctoqs}
Let $\bG=(V_\bG,E_\bG)$ be a graph that admits a good embedding with  length function $\ell$, and  $\on{carr}(\bG)= \bR^2$. Let $\sX$ denote the cable system of $\bG$. If $(\sX,d_1,\mu_1,\sE,\sF^1)$ is annular quasi-convex at large scales and satisfies \hyperlink{cap}{$(\on{cap}_{\le})$}, then $d_1$ and $d_\ell$ are quasisymmetric.
\end{theorem}
%
%
	We recall the definition of a weak quasisymmetry.
	\begin{definition} \label{d:wqs}
		{\rm 
			Given a homeomorphism $f:(\sX_1,d_1) \to (\sX_2,d_2)$,  $x \in \sX_1$ and $r>0$, set
			\be  \label{e:defH}
			H_f(x,r) = \frac{\sup\set{d_2(f(x),f(y)): d_1(x,y) \le r}}{\inf \set{d_2(f(x),f(y)):d_1(x,y) \ge r}}.
			\ee
			We say that $f:(\sX_1,d_1) \to (\sX_2,d_2)$ is a \emph{weak quasisymmetry}, if it is a   homeomorphism and if  there exists a 
			$H < \infty$ so that
			\be  \label{e:wqs}
			H_f(x,r) \le H
			\ee
			for all $x \in \sX_1, r>0$.
		}\end{definition}
		Every quasisymmetry is a weak quasisymmetry but the converse is not true in general \cite[Exercise 10.5]{Hei}. 
		Nevertheless, every weak quasisymmetry between geodesic metric spaces is a quasisymmetry \cite[Theorem 6.6]{Vai}.

{\sm {\em Proof of Theorem \ref{t:qctoqs}.}}
We follow the approach of \cite[Theorem 4.7]{HK} at large scales but we use a different argument using the goodness of embedding at small scales.

Let $f_\ell:(\sX,d_\ell) \to (\sX,d_1)$ denote the identity homeomorphism. 
By \cite[Theorem 6.6]{Vai}, it suffices to show that $f_\ell$ is a weak quasisymmetry.
Since the embedding is good, for any $R>0$, there exists $C_R \in (1,\infty)$ such that
for all $x,y \in \sX$ with $d_1(x,y) \le R$, we have
\be \label{e:qcs1}
C_R^{-1} r_x d_1(x,y) \le d_\ell(x,y) \le C_R r_x d_1(x,y),
\ee
where $r_x$ is the radius of separation.

Weak quasisymmetry of $f_\ell$ is equivalent to the following statement: there exists $H>0$
\be \label{e:qcs2}
d_\ell(x,a) \le d_\ell(x,b) \hspace{2mm} \mbox{ implies } \hspace{2mm} d_1(x,a) \le H d_1(x,b)  \hspace{2mm} \mbox{ for all $x ,a ,b \in \sX$.}
\ee
The estimate \eqref{e:qcs2} easily follows from \eqref{e:qcs1} in the case $d_1(x,a) \le R$ for any fixed $R$. Therefore it suffices to consider the case $d_1(x,a) \ge R$ for some large $R \in (1,\infty)$.

Suppose that
\be \label{e:qcs3}
 d_1(x,a) \ge R,\hspace{1mm} d_\ell(x,a) \le d_\ell(x,b), \hspace{1mm} \mbox{ and }  \hspace{1mm} s= d_1(x,a) > M d_1(x,b),
\ee
where $R,M$ will be chosen below. We will show that $M$  cannot be too large.

Choose a point $z \in \sX$ such that 
\bes 
d_\ell(x,z) \ge 2 d_\ell(x,b) \hspace{2mm} \mbox{and} \hspace{2mm} d_1(x,z) \ge s.
\ees
To see the existence of such a point $z$, we note that $\ol{B_\ell(x, 2 d_\ell(x,b))} \cup \ol{B_1(x,s)}$ is a compact subset of a non-compact space $\sX$, and therefore is a proper subset.

Let $E \subset B_1(x,s/M)$ denote the image of a shortest path in the $d_1$ metric joining $x$ to $b$.

 Since $(\sX,d_1)$ is annular quasi-convex at large scales, we can join $a$ and $z$ by a curve whose image $F$ is contained in $\sX \setminus B_1(x,s/C_L)$ provided $s \ge R_L$, where $C_L,R_L$ are constants from Definition \ref{d:llc}.

By \eqref{e:bpm4}, Lemma \ref{l:series}, there exists $C_1,K_1>0$ such that
\be \label{e:qcs4}
\Mod(E,F) \le \Mod( B_1(x,s/M), \sX \setminus B_1(x,s/C_L)) \le C_1 \left( \log (M/C_L) \right)^{-1},
\ee
provided $R \ge R_L \vee M$ and $M \ge K_1 C_L$.

Furthermore, $\dist_1(E,F) \ge s( C_L^{-1}- M^{-1}) \ge R( C_L^{-1}- M^{-1})$.
Note that, in the $d_\ell$ metric 
\[
\frac{\on{dist}_\ell(E,F)}{  \diam_{\ell}(E) \wedge \diam_\ell(F) } \le \frac{d_\ell(x,a)}{ d_\ell(x,b)} \le 1.
\]
By the large scale Loewner property of $(\sX,d_\ell,\mu_\ell)$, there exists $\epsilon,C_2 >0$ such that
\be \label{e:qcs5}
\Mod(E,F) \ge \delta
\ee
provided $R( C_L^{-1}- M^{-1})  \ge C_2$.
To arrive at a contradiction to \eqref{e:qcs3} from \eqref{e:qcs4} and \eqref{e:qcs5}, we require
\[
 \hspace{1mm} \log(M/C_L) \ge 2 C_1 \delta^{-1}, \hspace{1mm}  M \ge K_1 C_L, \hspace{1mm} \mbox{ and }  R \ge M \vee R_L \vee \left(  C_2 \left( C_L^{-1}- M^{-1}\right)^{-1} \right).
\]
The above requirements are clearly feasible by choosing $M$ using the first two constraints and then finally choosing $R$ using the third constraint above.
The contradiction to \eqref{e:qcs3}, along with \eqref{e:qcs1} implies that there are constants $C,R>0$ such that
\eqref{e:qcs2} holds with $H= M \vee (C_R^2)$.

\qed
\section{Heat kernel bounds using quasisymmetry} \label{s:hk}


In this section, we show  a characterization of the heat kernel estimates $\on{HK(\Psi)}$ corresponding to the cable process for a suitable $\Psi$ -- see Definition \ref{d:hkpsi}.

Cable systems $(\sX,d_1,\mu_1,\sE,\sF^1)$ corresponding to  a bounded degree graph  satisfies good local regularity properties that we summarize below.
\begin{lemma} (See \cite[Lemma 4.22]{BM2}(b)) \label{l:locreg}
Let  $(\sX,d_1,\mu_1,\sE,\sF^1)$ denote the cable system corresponding to a bounded degree graph $\bG$. Then for all $R>0$, there exists $C$ which only depends on $R$ and the bound on the degree, such that for all $x \in \sX, r \in (0,R]$ and $u \in \sF$,
\begin{align*}
C^{-1} r \le V_1(x,r) &\le Cr \\
\inf_{\alpha \in \bR} \int_{B_1(x,r)} \abs{u - \alpha}^2 \, d\mu_1 &\le C r^2 \int_{B_1(x,r)} \abs{\nabla_1 u}^2 \, d\mu_1, \\
C^{-1} \frac{1}{r} \le \Cap_{B_1(x,2r)}(B_1(x,r)) &\le C \frac{1}{r}.
\end{align*}
\end{lemma}
In other words, we have volume doubling property, Poincar\'e inequality and Capacity estimates on annuli for small balls by Lemma \ref{l:locreg}.

The following  theorem is the main result of this work.
\begin{theorem}  \label{t:main1}
Let $\bG$ be the planar graph of polynomial volume growth with volume growth exponent $d$.
Let $(\sX,d_\ell,\mu_\ell,\sE,\sF_\ell)$ be the MMD space corresponding to a good embedding of $\bG$ with carrier $\bR^2$. 
Then, the following are equivalent:
\begin{itemize}
\item [(a)] The metrics $d_1$ and $d_\ell$ are quasisymmetric.
\item [(b)] $(\sX,d_1,\mu_1,\sE,\sF^1)$ satisfies sub-Gaussian heat kernel bounds $\on{HKE}(\Psi)$ with $\Psi(r)= r^2 \vee r^d$.
\item [(c)]   $(\sX,d_1,\mu_1,\sE,\sF^1)$ satisfies  Poincar\'e inequality \hyperlink{pi}{$\on{PI}(\Psi)$} with $\Psi(r)=r^2 \vee r^d$, and the capacity bound \hyperlink{cap}{$(\on{cap}_\le)$}.
\item [(d)] $(\sX,d_1,\mu_1,\sE,\sF^1)$ is annular quasi-convex at large scales and   satisfies \hyperlink{cap}{$(\on{cap}_\le)$}.
\end{itemize}
\end{theorem}

\begin{remark} \label{r:main}
	{\rm
		\begin{enumerate}[(i)]
			\item  Theorem \ref{t:main1} can be generalized by replacing the assumption of polynomial volume growth with the volume doubling property (with respect to the counting measure).
			The space-time scaling function $\Psi$ in (b) and (c) above should be replaced by
			$\Psi(x,r)=r^2 \vee  V_\bG(x,r)$, where $V_\bG(x,r)$ denote the cardinality of $B_\bG(x,r)$. The proof of Theorem \ref{t:main1} easily extends to the general case using the methods in \cite[Section 5]{BM1}.
			\item 
The implication (c) $\Rightarrow$ (b) in Theorem \ref{t:main1} is essentially same as the conjecture in \cite[p. 1493]{GHL}, which we verify in a restricted setting. This conjecture has come to be known as the \emph{resistance conjecture}.
When the spectral dimension $d_s <2$ such an implication follows from the main results in \cite{BCK}. However, when the spectral dimension $d_s \ge 2$ the existing characterizations of heat kernel estimates seem too difficult to verify in practice -- see  the ICM survey \cite[Open problem III]{Kum} for further details. An alternate proof of the implication (c) $\Rightarrow$ (b) in a more general setting $d_f<1+d_w$ without the restriction of being planar is obtained in \cite{Mur}.
		\end{enumerate}
}\end{remark}
{\sm {\em Proof of Theorem \ref{t:main1}.}} Since the carrier is $\bR^2$, $(\sX,d_\ell)$ is complete. We will implicitly use the volume growth estimate of the cable system $(\sX,d_1,\mu_1)$ from \eqref{e:llc0}. \\
(a) $\Rightarrow$ (b):
By  \cite[Theorem 2.15]{BBK}, it suffices to verify the elliptic Harnack inequality and two sided bounds on the resistance of annuli. 

 By Theorem \ref{t:abgn} and \cite[Theorem 3.5]{St}, we obtain the parabolic Harnack inequality $\on{PHI}(2)$ for $(\sX,d_\ell,\mu_\ell,\sE,\sF^\ell)$. The parabolic Harnack inequality implies the elliptic Harnack inequality (EHI) for $(\sX,d_\ell,\mu_\ell,\sE,\sF^\ell)$ \cite[Propostion 3.2]{St}. 
As mentioned in Remark \ref{r:timechange}, $(\sX,d_\ell,\mu_\ell,\sE,\sF^\ell)$ is a time change of $(\sX,d_1,\mu_1,\sE,\sF^1)$ and therefore has the same harmonic functions.
By the quasisymmetry invariance of EHI, we obtain EHI for the time changed process $(\sX,d_1,\mu_1,\sE,\sF^1)$ \cite[Lemma 5.3]{BM1}.

The two sided bounds on the resistance of annuli in $(\sX,d_\ell,\mu_\ell,\sE,\sF^\ell)$ follows from 
 $\on{PHI}(2)$ and \cite[Theorem 2.15]{BBK}. 
 In particular, we have 
\[
\Cap_{B_\ell(x,2r)} (B_\ell(x,r)) \asymp 1, \hspace{3mm} \mbox{ for all $x \in \sX, r \ge r_x$.}
\] 
Using quasisymmetry, we can transfer the capacity bounds on annuli from one space to the other.
 By the same argument as  \cite[(5.20),(5.21), proof of the Theorem 5.14]{BM1} (see also \cite[Lemma 4.18(d)]{BM2}), we obtain 
\[
\Cap_{B_1(x,2r)} (B_1(x,r)) \asymp 1, \hspace{3mm} \mbox{ for all $x \in \sX, r \ge 1$.}
\] 
Therefore by \cite[Theorem 2.15]{BBK} and Lemma \ref{l:locreg}, we have (b). \\
(b) $\Rightarrow$ (c): This is immediate from \cite[Theorem 2.15 and 2.16]{BBK}.\\
(c) $\Rightarrow$ (d): This follows from Proposition \ref{p:llc}.\\
(d) $\Rightarrow$ (a): This is Theorem \ref{t:qctoqs}. \qed

{\sm {\em Proof of Theorem \ref{t:main2}.}} 
 By the results of \cite{BB}, the sub-Gaussian estimates  for the cable system $(\sX,d_1,\mu_1,\sE,\sF_1)$ are equivalent to sub-Gaussian estimates for the simple random walk.
The result now follows from the ring lemma \cite{RS}, comparison between the Euclidean metric and $d_\ell$ in \cite[Proposition 2.5]{ABGN}\footnote{Strictly speaking, \cite[Proposition 2.5]{ABGN} is stated for the case $\on{carr}(\bG) =\bU$ but the proof there works for the case $\on{carr}(\bG)=\bR^2$ as well.}, and the equivalence between (a) and (b) in Theorem \ref{t:main1}. 
\qed
\section{Examples} \label{s:xm}
\begin{figure}[h] 
	\centering
	\includegraphics[width=0.8\textwidth]{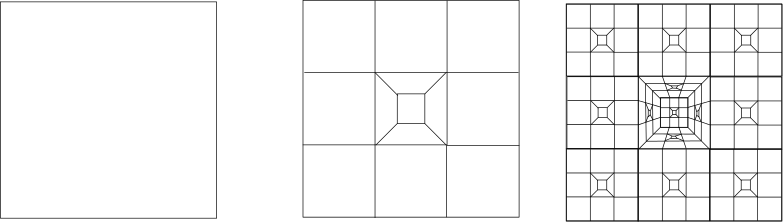}
	\caption{The sequence of graphs viewed from the central square converges to an infinite quadrangulation of the plane with walk dimension $d_w=\log_3(13)$.}
	\label{f:13-3}
\end{figure}

As an application of Theorem \ref{t:main1}, we present a new family of graphs that satisfy sub-Gaussian estimates. 
These graphs can be viewed as discrete analogues of fractal surfaces which we next describe.
\begin{example}[Snowball] \label{x:sb}
{\rm
Snowballs are fractals that are homeomorphic to $\bS^2$ and are defined as limits of polyhedral complexes. 
Their name stems from the fact that snowballs can be viewed as  higher dimensional analogues of the Koch snowflake. 
We recall the definition of one such fractal below.

Let $(\sS_0,d_0)$ denote the surface of the unit cube, equipped with the intrinsic metric $d_0$. In other words, $(\sS_0,d_0)$ can be viewed a polyhedral complex obtained by gluing six unit squares similar to the faces of a cube \cite[Defintion 3.2.4]{BBI}. We replace each \emph{face} in $\sS_0$ by 13 squares with edge length $\frac{1}{3}$ as shown in Figure \ref{f:13-3} to obtain a polyhedral complex $(\sS_1,d_1)$. More generally, we repeat this construction to obtain a geodesic metric space $(\sS_{n},d_{n})$ from $(\sS_{n-1},d_{n-1})$ ($n \ge 1$) by replacing each square of length $3^{-(n-1)}$ with 13 squares of each with edge length $3^{-n}$ as shown in Figure \ref{f:13-3}.
The polyhedral complex $(\sS_n, d_n)$ is obtained by gluing $6 \times (13)^n$ faces, where each face is isometric to a square with edge length $3^{-n}$.
It is easy to see that the metric spaces $(\sS_n,d_n)$ has a Gromov-Hausdorff limit $(\sS,d_\sS)$, which is called the snowball.
}\end{example}

We collect some properties of the metric space $(\sS,d_\sS)$.
Evidently, the spaces $(\sS_n,d_n)$ for $n \ge 0$, and $(\sS,d_\sS)$ are all homeomorphic to $\bS^2$. Let $d_{\bS^2}$ denote the standard Riemannian metric on $\bS^2$, viewed as an embedded surface in $\bR^3$.
It is known that $(\sS,d_\sS)$ and $(\bS^2,d_{\bS^2})$ are quasisymmetric -- see \cite{Mey02,Mey10} and \cite[p. 181]{BK}. 
We recall the following result of D.~Meyer that is essentially contained in \cite{Mey10}. 
\begin{prop} \label{p:mey}
There exists a homeomorphism $\eta:[0,\infty) \to [0,\infty)$, and $\eta$-quasisymmetric homeomorphisms $f_n: (\sS_n, d_n)\to (\bS^2,d_{\bS^2})$, $n \ge 0$, and $f:(\sS,d) \to (\bS^2,d_{\bS^2})$ satisfying the following properties: 
\begin{enumerate}[(a)]
	\item 
The push-forward metrics $\rho_n: \bS^2 \times \bS^2 \to [0,\infty)$, $n \ge 0$, where
$\rho_n(x,y)= d_n(f_n^{-1}(x),f_n^{-1}(y))$,
converge uniformly in $\bS^2 \times \bS^2$ to $\rho: \bS^2 \times \bS^2 \to [0,\infty)$, where $\rho(x,y)=d_{\sS}(f^{-1}(x),f^{-1}(y))$, that is 
\[
\lim_{n \to \infty}\sup_{x,y \in \bS^2} \abs{\rho_n(x,y) - \rho(x,y)}=0.
\]
\item Let $F^n_i, i=1,2,\ldots, 6 \cdot(13)^n$ denote the faces of the polyhedral complex $\sS_n, n \ge 0$.  We have
\be \label{e:mey}
\lim_{n\to \infty} \max_{1\le i \le 6(13)^n} \diam (f(F_i^n)) =0,
\ee 
where $\diam$ above denotes the diameter in $d_{\bS^2}$ metric. For $n \ge 1$, and for
$i=1,2,\ldots,  6 \cdot(13)^n$, there exists $j=1,\ldots, 6  \cdot(13)^{n-1}$ such that
\[
f_n(F^n_i) \subset f_{n-1}(F^{n-1}_j).
\]
\item The maps $f_n:\sS_n \to \bS^2$ are conformal maps when $\sS_n$ and $\bS^2$ are viewed as Riemann surfaces (The polyhedral surface $\sS_n$ has a canonical Riemann surface structure as explained in \cite[Section 3.3]{Bea}).
\end{enumerate} 
\end{prop}
	Next, we define a graph analogue  of  the snowball $\sS$ in Example \ref{x:sb}.
\begin{example}[Graphical snowball] \label{x:gsb}
	{\rm 
 We define the graph as a limit of finite graphs. We first define a finite planar graph $G_n=(V_n,E_n), n \ge 0$ using the polyhedral complexes $(\sS_n,d_n)$ defined in Example \ref{x:sb}.  
 The vertex set $V_n$ is same as the vertices of the polyhedron $\sS_n$, and two vertices $u,v \in V_n$ form an edge if and only if $d_n(u,v)= 3^{-n}$.
 Let $d^G_n$ denote the combinatorial graph metric on $V_n$. Let $p_n \in V_n$ be an arbitrary vertex in one of the six \emph{central faces} (there are 24 such vertices) -- see Figure \ref{f:13-3}. Then the sequence of pointed metric spaces $(V_n,d^G_n,p_n), n \ge 0$ has a pointed Gromov-Hausdorff limit as $n \to \infty$, $(V,d^G,p)$, where the metric $d^G$ can be viewed as the graph distance on a one-ended planar graph $G=(V,E)$ with volume growth exponent $d_f= \log_3(13)$.
 We call the graph $G=(V,E)$ the graphical snowball. 
 
 We claim that $G$ satisfies the equivalent conditions (a)-(d) in Theorem \ref{t:main1}.
Next, we sketch the proof of property (d) in Theorem \ref{t:main1}: annular quasi-convexity at large scales and the capacity upper bound \hyperlink{cap}{$(\on{cap}_\le)$}. 
  The  annular quasi-convexity at large scales easily follows from the corresponding
   property of the snowball $(\sS,d_\sS)$. The proof of \cite[Proposition 18.5(ii)]{BK}
    can be easily adapted to the graph setting using the comparison between the intrinsic metric $d_\sS$ and the `visual metric' in \cite[Lemma 2.2]{Mey02}.

  The estimate on capacity \hyperlink{cap}{$(\on{cap}_\le)$} for the graphical snowball $G$ is obtained using modulus estimates and comparison of modulus between metric spaces and their discrete graph approximations in \cite{BK}.  Next, we sketch the proof of the capacity upper bound \hyperlink{cap}{$(\on{cap}_\le)$}.
  Let $\wt{G}_n$ denote the face barycenter triangulation of $G_n$ (see \cite{CFP2} for the definition of face barycenter triangulation). Then $\wt{G}_n$ is a $K$-approximation on $(\sS,d_\sS)$ in the sense of {BK}.  Roughly speaking, a $K$-approximation is a covering of the space indexed by the vertices of a graph, such that the covering has controlled overlap, that each set in the covering is approximately a ball, and adjacent vertices have correspond to sets that intersect with comparable sizes (see \cite[p. 141]{BK} for precise definition). In our case, we can choose $\wt{G}_n$ to cover $\sS$ by balls of radius $3^{-n+2}$.
 By  \cite[Theorem 11.1]{BK}, we  obtain (combinatorial) modulus estimates on the annuli of $\wt{G}$, where $\wt{G}$ denotes the face barycenter triangulation of the graphical snowball $G$.
  There are two different notions of (combinatorial) modulus in the context of graphs, one of which assigns weights to edges and the other assigns weights to vertices -- see the definition of vertex extremal length and edge extremal length  in \cite[p. 128]{HS} (extremal length is the reciprocal of  modulus). The notion of modulus used in \cite{BK} assigns weights to vertices. 
    However, for the capacity bounds, the version of modulus that assigns weights to edges is relevant \cite[p. 128]{HS}. 
    For bounded degree graphs, the two versions of modulus are comparable up to a multiplicative factor (that depends only the uniform bound on the degree) -- \cite[proof of Theorem 8.1]{HS}. Furthermore, since $\wt{G}$ and $G$ are quasi-isometric graphs, the modulus of annuli are comparable.
    Combining the above observations, we obtain \hyperlink{cap}{$(\on{cap}_\le)$} in Theorem \ref{t:main1}(d). Hence, we obtain the following: 
    \begin{prop} \label{p:gsb}
     	The graphical snowball $G$ has polynomial  growth with volume growth exponent $d=\log_3(13)$ and satisfies sub-Gaussian heat kernel bounds with $d_w=\log_3(13)$.
    \end{prop}
  
  We remark that the choice of base points $p_n, n \ge0$ made in the definition of $G$ is for  concreteness, and the properties we discussed above is independent of this choice. Since the graphs $G_n$ have uniformly bounded degree, any such sequence of pointed metric spaces will have a sub-sequential limit, that can be viewed as an infinite graph. 
}\end{example}

The snowball and its graph version presented in Examples \ref{x:sb} and \ref{x:gsb} admit many  variants, which also have the same quasisymmetry property; see \cite[Theorem 1A, Remarks on p. 1268]{Mey10} for details.
The graph mentioned in Example \ref{x:gsb}  can be viewed as a  net of a metric tangent cone (see \cite[Definition 8.2.2]{BBI}). Yet another viewpoint is that these are graph versions of expansion complexes corresponding to a finite subdivision rule \cite{CFP1,CFP2} -- see also \cite{BS,CFKP}. 
\begin{example}[`Regular' Pentagonal tiling]{\rm
	The following example is a graph version of a Riemann surface studied in
 \cite{BS,CFP1,CFKP}. 
 The fractal analogue of Example \ref{x:sb} is built using a dodecahedron where each of the twelve pentagonal faces is subdivided into six pentagons as shown in Figure \ref{f:6-2}, where at the $n$-th iteration, the polyhedral surface is obtained by gluing 
 $12\times(6)^n$ pentagons where the length of each side is $2^{-n}$ (here $n=0$ corresponds to the dodecahedron). 
 The graph analogue of Example \ref{x:gsb} can be obtained using a similar construction where the base point is chosen from one of the vertices in the central pentagon.  The resulting graph satisfies sub-Gaussian heat kernel estimates with spectral dimension $2$ and  with  walk dimension $d_w=\log_2(6)$ using the same methods discussed in Examples \ref{x:gsb} and \ref{x:sb}. 
 
 We remark that the question of finding the spectral dimension of this example was raised by B\'alint T\'oth in 1988 \cite{Kum}, \cite{Tot}.  We refer the reader to \cite[Section 6]{Mey02} for a zoo of examples. For instance, the Xmas tree in \cite[Section 6.4]{Mey02} has qualitative resemblance to fractals found in nature such as broccoli or cauliflower.
	 \begin{figure}[h] 
	 	\centering
	 	\includegraphics[width=0.8\textwidth]{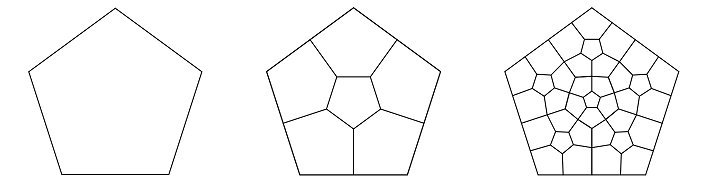}
	 	\caption{The sequence of graphs viewed from the central pentagon converges to an infinite pentagulation of the plane with  walk dimension $d_w=\log_26$.}
	 	\label{f:6-2}
	 \end{figure}
}\end{example}

\subsection{Diffusion on the snowball}
We show that the ideas presented in the earlier sections for random walks on graphs also apply to diffusions on fractals. We explain this on the snowball $(\sS,d_{\sS})$ defined in Example \ref{x:sb}. We will define the canonical diffusion on the snowball as a limit of diffusions on its polyhedral approximations $\sS_n, n \ge 0$ as $n \to \infty$.

We need a measure $\mu$ on $\sS$ that plays the role of a symmetric measure for the diffusion. We again construct $\mu$ as a limit of measures defined on polyhedral approximations.
There is a natural family of measures on $\sS_n, n \ge0$ and $\sS$  that we next describe.
Let $\mu_n$ denote the surface area measure on $(\sS_n,d_n)$ normalized to be a probability measure (or equivalently, the normalized Riemannian measure on $\bS^2$), so that each face has rescaled Lebesgue measure with  $\mu_n(F_i^n)={13^{-n}}/{6}$ for all $i=1,\ldots,6\cdot(13)^n$. Let $f_n,f$ denote the quasisymmetric maps in Proposition \ref{p:mey}.
The push-forward measures $\nu_n= (f_n)_* \mu_n, n \ge 0$ on $\bS^2$, converge weakly as $n \to \infty$ to a measure $\nu=f_* \mu$ where $\mu$ is a Borel probability measure on $(\sS,d_\sS)$. The weak convergence is an easy consequence of Proposition \ref{p:mey}(ii).

We state some elementary estimates on the measures $\mu_n, \mu$ defined above.
Let $B_n$ and $B_\sS$ denote the open balls in $(\sS_n,d_n)$ and $(\sS,d_\sS)$  respectively and let $\phi_n:[0,\infty) \to [0,\infty)$ denote the function
\[
\phi_n(r)=\begin{cases}
\frac{9^n}{13^n} r^2 & \mbox{if $r \le 3^{-n}$},\\
r^{\log_3(13)} &\mbox{if $3^{-n}\le r \le 1$,}\\
1 & \mbox{if $r \ge 1$.}
\end{cases}
\]
Furthermore, 
the measure $\mu$ is $\log_3(13)$-Ahlfors regular; that is, there exists $C>1$ such that for all $0<r \le \diam(\sS,d_\sS)$ and for all $x \in \sS$, we have
\be \label{e:reg1}
C^{-1} r^{\log_3(13)} \le  \mu\left(B_{\sS}(x,r)\right) \le C r^{\log_3(13)},
\ee 
and for all $n \ge 0$, for all  $0<r \le \diam(\sS_n,d_n)$ and for all $x \in \sS_n$, we have
\be \label{e:reg2}
C^{-1} \phi_n(r) \le \mu_n(B_n(x,r)) \le C \phi_n(r).
\ee 
The above volume estimates follows from comparing the balls $B_{\sS}(x,r)$ to faces of $F^n_i$ where $n$ is chosen such that $3^{-n} \approx r$ and using \cite[Lemma 2.2]{Mey02}. For instance, the proof of \cite[Proposition 18.2]{BMe} can be easily adapted to show \eqref{e:reg1} and \eqref{e:reg2}.

Recall that $\nu_n=(f_n)_* \mu_n, \nu=f_*\mu$ are measures on $\bS^2$, where $f_n:\sS_n \to \bS^2, f: \sS \to \bS^2$ denote the maps defined  in Proposition \ref{p:mey}. Let $(\sE^n,\sF^n)$ denote the Dirichlet form corresponding to the Brownian motion on the Euclidean complex $(\sS_n,d_n)$ equipped with the symmetric measure $\mu_n$.
The estimates \eqref{e:reg1}, \eqref{e:reg2}, along with Proposition \ref{p:mey}  imply that the family of metric measure spaces $(\bS^2,d_{\bS^2},\nu_n), n \ge 0$, and  $(\bS^2,d_{\bS^2},\nu)$ are \emph{uniformly doubling}, that is there exists $C_D>1$ such that
\be \label{e:unifd}
\nu_n(B_{\bS^2}(x,2r)) \le C_D \nu_n(B_{\bS^2}(x,r)), \mbox{ and }\nu(B_{\bS^2}(x,2r)) \le C_D \nu(B_{\bS^2}(x,r)), 
\ee
for all $x \in \bS^2, r >0$.
In \cite[Lemma 1.14 and 1.15]{PS}, the authors provide two equivalent definitions of  $(\sE^n,\sF^n)$. 
Let $(X^n_t)_{t \ge 0}$ denote the process corresponding to the Dirichlet form $(\sE^n,\sF^n)$ on $L^2(\sS_n,\mu_n)$. Let $(Y^n_t)_{t \ge 0}$ denote the process on $\bS^2$ corresponding to the image of $X_t^n$ under the homeomorphism $f_n$, that is $Y^n_t=f_n(X^n_t)$ for all $t \ge 0$. By the conformal invariance of Brownian motion, we have that the Dirichlet form $(\wt{\sE}^n,\wt{F}^n)$ on $L^2(\bS^2,\nu_n)$ corresponding the process  $(Y^n_t)_{t \ge 0}$ is 
\be \label{e:cinv}
\wt{\sE}^n(g,g)= \int_{\bS^2} \abs{\nabla g}^2 \, dm, \hspace{3mm} \mbox{for all $g \in \wt{\sF}^n= \{ h \circ f_n^{-1}: h \in \sF^n\}$},
\ee 
where $\abs{\nabla g}$, and $m$ denotes the length of the  Riemannian gradient, and Riemannian measure respectively on $\bS^2$. In other words, $(Y^n_t)_{t \ge 0}$ is a time change of the Brownian motion on $\bS^2$ using the measure $\bS^2$. 
Although the measure $\mu_n$ is not $\mathcal{C}^\infty$ with respect to the  standard atlas on $\bS^2$, it is  $\mathcal{C}^\infty$ on $\bS^2$ except for finitely many points that are images of the vertices of $\sS_n$ under $f_n$. 
Since any finite set has capacity zero on $\bS^2$, the conformal invariance  in \eqref{e:cinv} follows even though the measure is not $\mathcal{C}^\infty$.
By \eqref{e:cinv}, the diffusions $(X^n_t)_{t \ge 0}$ on $\sS_n$ can be viewed as time change of Brownian motion on $\bS^2$ with Revuz measure $\nu_n$ using the map $f_n, n \ge 0$. 
By Proposition \ref{p:mey}, the measures $\nu_n$ converge weakly to $\nu$ on $\bS^2$ as $n \to \infty$. 
\begin{lemma}  \label{l:revuz}
Let $(\wt{\sE},\wt{\sF})$ on $L^2(\bS^2,m)$, denote the Dirichlet form corresponding to the
Brownian motion on $\bS^2$. Let 
\[
\wt{\sE}_1(f,f)= \wt{\sE}(f,f) + \norm{f}_{L^2(m)}^2.
\]
Then the measure $\nu=f_* \mu$ defined above is of finite energy integral: there exists $C>0$ such that
\be \label{e:fe}
\int_{\bS^2} \abs{g(x)}\,\nu(dx) \le C \sqrt{ \wt{\sE}_1(g,g)}, \hspace{3mm} \mbox{for all $g \in \wt{\sF}\cap \mathcal{C}_0(\bS^2)$}.
\ee
In particular, $\nu$ charges no set of zero capacity. Furthermore, $\nu$ has full support.
\end{lemma}
\proof
The fact that $\nu$ charges no set of zero capacity follows from \eqref{e:fe} and \cite[Lemma 2.2.3]{FOT}. Evidently, \eqref{e:reg1} and the quasisymmetry of $f:(\sS,d_\sS) \to (\bS^2,d_{\bS^2})$ imply that $\nu$ has full support.

It only remains to verify \eqref{e:fe}.
Let $p_t:\bS^2 \times \bS^2 \to [0,\infty]$ denote the continuous version of the heat kernel of  
the Brownian motion on $(\bS^2,d_{\bS^2},m)$ and let $G_1:\bS^2 \times \bS^2 \to [0,\infty]$ denote the massive Green function
\[
G_1(x,y)= \int_0^\infty p_t(x,y) e^{-t} \, dt.
\]
 Using Gaussian estimates for $p_t$, we obtain the following estimate: there exist $C_1>0$ such that
 \be \label{e:sm1}
 G_1(x,y) \le C_1 \left( \ln \frac{1}{d_{\bS^2}(x,y)} + C_1 \right).
 \ee
Define the $1$-potential of the measures $\nu_n$ and $\nu$ as
\[
g_1(\nu)(x):= \int_{\bS^2} G_1(x,y) \, \nu(dy), \hspace{4mm} g_1(\nu_n)(x):= \int_{\bS^2} G_1(x,y) \, \nu_n(dy),
\]
for all $n\ge 0$, and $x \in \bS^2$.
Since $\nu_n \ll m$ and $\frac{d \nu_n}{dm} \in L^\infty(m)$, by \cite[(1.3.1) and (2.2.2)]{FOT} we have that $g_1(\nu_n)$ is the $1$-potential $U_1 \nu_n, n \ge 0$, that is $g_1(\nu_n)=U_1 (\nu_n) \in \wt{\sF}$.

Next, we show that the functions $g_1(\nu_n), n \ge 0$ and  $g_1(\nu)$ admit continuous versions. By the uniform doubling property \eqref{e:unifd} of the measures $\nu$ and $\nu_n, n \ge 0$, along with the reverse volume doubling property in \cite[Exercise 13.1]{Hei}, there exists $C_2>1, \alpha >0$ such that
\be \label{e:sm2}
 \sup_{n \ge 0}  \sup_{x \in \bS^2} \sup_{ r \in (0, \diam(\bS^2,d_{\bS^2}))}  r^{-\alpha} \left( \nu_n(B_{\bS^2}(x,r)) + \nu(B_{\bS^2}(x,r)) \right) \le C_2
\ee
Using \eqref{e:sm1}, \eqref{e:sm2}, and the same argument  as in \cite[Proposition 2.3]{GRV}, 
we obtain that for all $n \ge 0$, $g_1(\nu_n) \in \mathcal{C}(\bS^2) \cap \wt{\sF}, g_1(\nu) \in \mathcal{C}(\bS^2)$ and that they are uniformly bounded; that is, there exists $C_3>0$ such that
\be \label{e:sm3}
 \sup_{n \ge 0} \sup_{x \ge 0} \left( g_1(\nu_n)(x)+ g_1(\nu)(x) \right) \le C_3.
\ee 
Furthermore using \eqref{e:sm1}, \eqref{e:sm2}, and a straightforward adaptation of the argument in \cite[Proposition 2.3]{GRV}, we have that $g_1(\nu_n)$ converges uniformly to $g_1(\nu)$ as $n \to \infty$, \emph{i.e.}, 
\be \label{e:sm4}
\lim_{n \to \infty} \norm{g_1(\nu_n)-g_1(\nu)}_\infty= \lim_{n \to \infty} \sup_{x \in \bS^2} \abs{g_1(\nu_n)(x) - g_1(\nu)(x)}  = 0.
\ee
By \cite[(2.2.2)]{FOT} and \eqref{e:sm3}, we obtain
\be \label{e:sm5}
\wt{\sE}_1(g_1(\nu_n),g_1(\nu_n)) \le \int_{\bS^2}  g_1(\nu_n)(x) \, \nu_n(dx)\le \nu_n(\bS^2) \norm{g_1(\nu_n)}_\infty \le C_3.
\ee
By \eqref{e:sm5}, \cite[(2.2.2)]{FOT}, and the Cauchy-Schwarz inequality, we obtain  for all $v \in \wt{\sF} \cap \sC_0(\bS^2)$, 
\[
\int_{\bS^2} \abs{v(x)}\, \nu(dx)= \lim_{n \to \infty} \int_{\bS^2} \abs{v(x)}\, \nu_n(dx) = \lim_{n \to \infty}\wt{\sE}_1(g_1(\nu_n),v)  \le \sqrt{C_3} \sqrt{\wt{\sE}_1(v,v)}.
\]
Therefore, $\nu$ is of finite energy integral. Using Lemma \ref{l:extd} and \eqref{e:sm5}, we obtain that $g_1(\nu) \in \sF \cap \sC_0(\bS^2)$ and that $g_1(\nu)$ coincides with the $1$-potential of the measure $\nu$.
\qed

By Lemma \ref{l:revuz} and the Revuz correspondence \cite[Theorem 5.1.4]{FOT}, there is a time change of the Brownian motion on $\bS^2$ that is symmetric with respect to $\nu$, which we denote by $(Y_t)_{t \ge 0}$. The process $(Y_t)$ can be viewed as the limit of $(Y^n_t)_{t \ge 0}$, as the measures $\nu_n$ converge to $\nu$ on $\bS^2$. Since $Y^n_t=f_n(X^n_t)$ is  the image of the diffusion on the polyhedral approximations $(\sS_n,d_n,\mu_n)$ of the snowball $(\sS,d_\sS,\mu)$, we can view $Y_t$ as the image of the canonical $\mu$-symmetric diffusion  $(X_t)_{t \ge 0}$ on $(\sS,d_\sS,\mu)$, where $X_t= f^{-1}(Y_t)$. This defines the diffusion on the snowball as the limit of diffusions on the polyhedral approximations $\sS_n$ as $n \to \infty$.

The proof of the implication `(a) $\Rightarrow$ (b)' in Theorem \ref{t:main1} can be readily adapted to obtain heat kernel bounds for the process $(X_t)_{t \ge 0}$ on $(\sS,d_{\sS}, \mu)$. Let $p^\sS_t:\sS \times \sS \to [0,\infty]$ denote the 
heat kernel for $(X_t)_{t \ge 0}$.
The Gaussian estimates for the Brownian motion on $\bS^2$ along with quasisymmetry of $f$, yields the following sub-Gaussian estimate on $p_t^{\sS} (\cdot, \cdot)$: there exists $C>1$ such that
\[
\frac{C^{-1}}{1 \vee t}  \exp \left( - \frac{ Cd_{\sS}(x,y)^{d/(d-1)}}{ t^{1/d}} \right) \le  p_t^{\sS}(x,y) \le \frac{C}{1 \vee t}  \exp \left( - \frac{ d_{\sS}(x,y)^{d/(d-1)}}{C t^{1/d}} \right)
\]
for all $t >0$ and for all $x,y \in \sS$, where $d=\log_3(13)$. In other words, the canonical diffusion on $(\sS,d_\sS,\mu)$ has spectral dimension $d_s=2$ and walk dimension $d_w=\log_3(13)$.

\med {\bf Acknowledgement.}
I am grateful to Martin Barlow for providing helpful suggestions on an earlier draft of this paper.
I thank Tim Jaschek for pointing out several typos in an earlier draft.
 I thank Steffen Rohde for his interest in this work,  encouragement, and several  discussions on the quasisymmetry of snowballs.
I  benefited from conversations with Omer Angel, and Tom Hutchcroft on circle packings, and with Jun Kigami on quasisymmetry. 
I thank Phil Bowers and Ken Stephenson for the permission to use Figure \ref{f:6-2} which appeared in \cite[p. 62]{BS}.
 I thank the anonymous referee for several helpful remarks and corrections.

\noindent Department of Mathematics, University of British Columbia,
Vancouver, BC V6T 1Z2, Canada. \\
mathav@math.ubc.ca

\end{document}